\documentclass[11pt, twoside, reqno, aap, preprint]{article}

\RequirePackage[OT1]{fontenc}
\RequirePackage[aos,amsthm,amsmath,noinfoline]{myimsart}
\RequirePackage[dvips,colorlinks]{hyperref}

\usepackage{amssymb}
\usepackage{amsopn}
\usepackage{mathrsfs}

\input epsf

\startlocaldefs
\theoremstyle{plain}
\newtheorem{theorem}{Theorem}[section]
\newtheorem{lemma}{Lemma}[section]
\newtheorem{proposition}{Proposition}[section]

\theoremstyle{definition}

\newtheorem{remark}{Remark}[section]
\newtheorem{corollary}{Corollary}[section]

\def\betam{\bar{\beta}}

\def\cc#1{\mathcal{#1}}
\newcommand{\B}{\cc{B}^{*}}
\def\T{\top}
\def\tr{\mathop{\rm tr}}
\def\Span{\mathop{\rm span}}
\def\rank{\mathop{\rm rank}}
\def\arginf{\mathop{\rm arg\,inf}}
\def\esp{{\bf E}}
\def\Pb{{\bf P}}

\def\RR{\mathbb R}

\def \1{{\rm 1}\mskip -4,5mu{\rm l} } 

\def\eps{\varepsilon}

\def\theta{\vartheta}

\def \1{{\rm 1}\mskip -4,5mu{\rm l}}
\newcommand{\onev}[1] {\biggl(\begin{matrix} 1 \\[-3pt] #1
\end{matrix}\biggr)}

\def\diag{\operatorname{diag}}

\newcommand{\M}{\mathscr M}
\newcommand{\tM}{\ \ \tilde{\!\!\!\!\!\mathscr M}}

\newcommand{\nablaf}{\nabla \! f}
\newcommand{\nablaG}{\nabla \! g}

\def\V{\tilde{V}}
\def\e{\boldsymbol{e}}

\newenvironment{Proof}{\noindent \bf Proof.\  \rm}{\hfill \qed}

\numberwithin{equation}{section}

\endlocaldefs

\begin{document}

\parskip5pt
\parindent0pt

\begin{frontmatter}

\title{A new algorithm for
estimating the effective dimension-reduction subspace}

\runtitle{Estimation of the dimension-reduction subspace}

\begin{aug}
\author{\fnms{Arnak} \snm{Dalalyan}\ead[label=e1]{dalalyan@ccr.jussieu.fr}},
\author{\fnms{Anatoly}
\snm{Juditsky}\ead[label=e2]{anatoli.iouditski@imag.fr}}
\and
\author{\fnms{Vladimir} \snm{Spokoiny}
\ead[label=e3]{spokoiny@wias-berlin.de}
}

\runauthor{Dalalyan, Juditsky and Spokoiny}

\affiliation{University Paris 6, University Joseph Fourier of Grenoble
and Weierstrass Institute for Applied Analysis and Stochastics}

\address{Laboratoire de Probabilit\'es\\
Universit\'e Paris 6, Bo\^\i te courrier 188\\
75252 Paris Cedex 05, France\\
\printead{e1}}

\address{LMC-IMAG\\
51 rue des Mathematiques,
B. P. 53\\
38041 Grenoble Cedex 9,
France\\
\printead{e2}}

\address{Weierstrass-Institute\\
Mohrenstr.\ 39\\
10117 Berlin Germany\\
\printead{e3}}

\end{aug}

\begin{abstract}
The statistical problem of estimating the effective
dimension-reduction (EDR) subspace in the multi-index regression
model with deterministic design and additive noise is considered. A
new procedure for recovering the directions of the EDR subspace is
proposed. Under mild assumptions, $\sqrt n$-consistency of the
proposed procedure is proved (up to a logarithmic factor) in the
case when the structural dimension is not larger than $4$. The
empirical behavior of the algorithm is studied through numerical
simulations.
\end{abstract}



\end{frontmatter}

\section{Introduction}\label{sec1}

One of the most challenging problems in modern statistics is to find
efficient methods for treating high-dimensional data sets. In
various practical situations the problem of predicting or explaining
a scalar response variable $Y$ by $d$ scalar predictors
$X^{(1)},\ldots, X^{(d)}$ arises. For solving this problem one
should first specify an appropriate mathematical model and then find
an algorithm for estimating that model based on the observed data.
In the absence of a priori information on the relationship between
$Y$ and $X=(X^{(1)},\ldots, X^{(d)})$, complex models are to be
preferred. Unfortunately, the accuracy of estimation is in general a
decreasing function of the model complexity. For example, in the
regression model with additive noise and two-times continuously
differentiable regression function $f :\RR^d\to\RR$, the most
accurate estimators of $f$ based on a sample of size $n$ have a
quadratic risk decreasing as $n^{-4/(4+d)}$ when $n$ becomes large.
This rate deteriorates very rapidly with increasing $d$ leading to
unsatisfactory accuracy of estimation for moderate sample sizes.
This phenomenon is called ``curse of dimensionality'', the latter
term being coined by Bellman (1961).

To overcome the ``curse of dimensionality'', additional restrictions
on the candidates $f$ for describing the relationship between $Y$
and $X$ are necessary. One popular approach is to consider the multi-index
model with $m^*$ indices: for some linearly independent vectors $\theta_1$,
$\ldots$, $\theta_{m^*}$ and for some function $g:\RR^{m^*}\to \RR$, the
relation
$f(x)=g(\theta^\top_1x,\ldots,\theta^\top_{m^*}x)$ holds for every
$x\in\RR^d$.
Here and in the sequel the vectors are understood as one column matrices
and $M^\top$ denotes the transpose of the matrix $M$.
Of course, such a restriction is useful only if $m^*<d$ and the main
argument
in favor of using the multi-index model is that for most data sets
the  underlying structural dimension $m^*$ is substantially smaller than $d$.
Therefore, if the vectors $\theta_1$, $\ldots$, $\theta_{m^*}$ are known,
the estimation of $f$ reduces to the estimation of $g$,
which can be performed much better because of lower dimensionality of the
function $g$ compared to that of $f$.

Another advantage of the multi-index model is that it assesses that
only few linear combinations of the predictors may suffice for
``explaining'' the response $Y$. Considering these combinations as
new predictors leads to a much simpler model (due to its low
dimensionality), which can be successfully analyzed by graphical
methods, see \cite{CookWeis}, \cite{Cook} for more details.

Thus, throughout this work we assume that we are given $n$ observations
$(Y_1,X_1),\ldots,(Y_n,X_n)$  from the model
\begin{equation}\label{model}
Y_i=f(X_i)+\varepsilon_i=g(\theta_1^{\T}X_i,\ldots,\theta_{m^*}^{\T}X_i)
+\varepsilon_i,
\end{equation}
where $\varepsilon_1,\ldots,\varepsilon_n$ are unobserved errors
assumed to be mutually independent zero mean random variables,
independent of the design $\{X_i,i\le n\}$.

Since it is unrealistic to assume that
$\theta_1,\ldots,\theta_{m^*}$ are known, estimation of these
vectors from the data is of high practical interest. When the
function $g$ is unspecified, only the linear subspace $\mathcal
S_{\theta}$ spanned by these vectors may be identified from the
sample. This subspace is usually called \textit{index space} or
\textit{dimension-reduction (DR) subspace}. Clearly, there are many
DR subspaces for a fixed model $f$. Even if $f$ is observed without
error, only the smallest DR subspace, henceforth denoted by
$\cc{S}$, can be consistently identified. This smallest DR subspace,
which is the intersection of all DR subspaces, is called
\textit{effective dimension-reduction} (EDR) subspace \cite{LiSIR}
or \textit{central mean subspace} \cite{CookLi}. We adopt in this
paper the former term, in order to be consistent with \cite{HJPS}
and \cite{XiaetAl}, which are the closest references to our work.

The present work is devoted to studying a new algorithm for
estimating the EDR subspace.We call it structural adaption via
maximum minimization (SAMM). It can be regarded as a branch of the
structure-adaptive (SA) approach introduced in \cite{HJS},
\cite{HJPS}.

Note that a closely related problem is the estimation of the central
subspace (CS), see \cite{CookWeis} for its definition. For model
(\ref{model}) with i.i.d.\ predictors, the CS coincides with the EDR
subspace. Hence, all the methods developed for estimating the CS can
potentially be applied in our set-up. We refer to \cite{CookLi} for
background on the difference between the CS and the central mean
subspace and to \cite{CookNi} for a discussion of the relationship
between different algorithms estimating these subspaces.

There are a number of methods providing an estimator of the EDR
subspace in our set-up. These include ordinary least square
\cite{LiDuan}, sliced inverse regression \cite{LiSIR}, sliced
inverse variance estimation \cite{CookWeis1}, principal Hessian
directions \cite{Li}, graphical regression \cite{Cook}, parametric
inverse regression \cite{BuraCook}, SA approach \cite{HJPS},
iterative Hessian transformation \cite{CookLi}, minimum average
variance estimation (MAVE) \cite{XiaetAl} and minimum discrepancy
approach \cite{CookNi}.

All these methods, except SA approach and MAVE,
are related to the principle of inverse regression (IR). Therefore
they inherit its well known limitations.
First, they require a hypothesis on the probabilistic
structure of the predictors usually called linearity condition. Second,
there is no theoretical justification guaranteeing that these methods
estimate the whole EDR subspace and not just a part thereof
(cf. \cite[Section 3.1]{CookLi04} and the comments on
the third example in \cite[Section 4]{HJPS}). In the same
time, they have the advantage of being simple for implementation
and for inference.

The two other methods mentioned above -- SA approach and MAVE --
have much wider applicability including even time series analysis.
The inference for these methods is more involved than that of IR
based methods, but SA approach and MAVE are recognized to provide
more accurate estimates of the EDR subspace.

These arguments, combined with the empirical experience, indicate
the complementarity of different methods designed to estimate the
EDR subspace. It turns out that there is no procedure among those cited
above that outperforms all the others in plausible settings.
Therefore, a reasonable strategy for estimating the EDR subspace is to
execute different procedures and to take a decision after comparing
the obtained results. In the case of strong contradictions,
collecting additional data or using extra-statistical arguments is
recommended.

The algorithm SAMM we introduce here exploits the fact that the
gradient $\nabla f$ of the regression function $f$ evaluated at any
point $x\in\RR^d$ belongs to the EDR subspace. The estimation of the
gradient being an ill-posed inverse problem, it is better to
estimate some linear combinations of $\nabla f(X_1),\ldots,\nabla
f(X_n)$, which still belong to the EDR subspace.

Let $L$ be a positive integer. The main idea leading to the
algorithm proposed in \cite{HJPS} is to iteratively estimate $L$
linear combinations $\beta_1,\ldots,\beta_L$ of vectors $\nabla
f(X_1),\ldots,\nabla f(X_n)$ and then to recover the EDR subspace from
the vectors $\beta_\ell$ by running a principal component analysis
(PCA). The resulting estimator is proved to be $\sqrt n$-consistent
provided that $L$ is chosen independently on the sample size
$n$. Unfortunately, if $L$ is small with respect to $n$, the subspace
spanned by the vectors $\beta_1,\ldots,\beta_L$ may cover only a
part of the EDR subspace. Therefore, the empirical experience advocates
for large values of $L$, even if the desirable feature of $\sqrt
n$-consistency fails in this case.

The estimator proposed in the present work is designed to provide a
remedy for this dissension between the theory and the empirical
experience. This goal is achieved by introducing a new method of
extracting the EDR subspace from the estimators of the vectors
$\beta_1,\ldots,\beta_L$. If we think of PCA as the solution to a
minimization problem involving a sum over $L$ terms (see (\ref{PCA})
in the next section) then, to some extent, our proposal is to
replace the sum by the maximum. This motivates the term structural
adaptation via maximum minimization. The main advantage of SAMM is
that it allows us to deal with the case when $L$ increases
polynomially in $n$ and yields an estimator of the EDR subspace
which is consistent under a very weak identifiability assumption. In
addition, SAMM provides a $\sqrt n$-consistent estimator (up to a
logarithmic factor) of the EDR subspace when $m^*\leq 4$.

If $m^*=1$, the corresponding model is referred to as
\textit{single-index} regression. There are many methods for
estimating the EDR subspace in this case (see \cite{YinCook},
\cite{DHP} and the references therein). Note also that the methods
for estimating the EDR subspace have often their counterparts in the
partially linear regression analysis, see for example \cite{SSV} and
\cite{ChLiT}.

Some aspects of the application of dimension reduction techniques in
bioinformatics are studied in~\cite{ALL} and \cite{BuPf}.

The rest of the paper is organized as follows. We review the
structure-adaptive approach and introduce the SAMM procedure in
Section~\ref{StrAd}. Theoretical features including $\sqrt
n$-consistency of the procedure are stated in Section~\ref{Theor}.
Section~\ref{Sim} contains an empirical study of the proposed
procedure through Monte Carlo simulations. The technical proofs are
deferred to the Appendix.

\section{Structural adaptation and SAMM}
\label{StrAd}

Introduced in \cite{HJS}, the structure-adaptive approach is based
on two observations. First, knowing the structural information helps
better estimate the model function. Second, improved model
estimation contributes to recovering more accurate structural
information about the model. These advocate for the following
iterative procedure. Start with the null structural information,
then iterate the above-mentioned two steps (estimation of the model
and extraction of the structure) several times improving the quality
of model estimation and increasing the accuracy of structural
information during the iteration.

\subsection{Purely nonparametric local linear estimation}
When no structural information is available,
one can only proceed in a fully nonparametric way.
A proper estimation method is based on local linear smoothing
(cf.\ \cite{FanGij} for more details):
estimators of the function
$f$ and its gradient $\nablaf$ at a point $X_{i}$
are given
by
\begin{align*}
\begin{pmatrix}
     \hat{f}(X_{i}) \\
     \widehat{\nablaf}(X_{i})
\end{pmatrix}
&=\arginf_{(a, c)^\T} \sum_{j=1}^{n} \big( Y_{j} - a - c^{\T} X_{ij}
  \big)^{2}K\big( |X_{ij}|^{2}/b^{2}\big)\nonumber\\
&=\bigg\{\sum_{j=1}^{n}\onev{X_{ij}}\onev{X_{ij}}^{\T}
    K\bigg(\frac{|X_{ij}|^{2}}{b^{2}}\bigg)\bigg\}^{-1}
    \sum_{j=1}^{n}Y_{j} \onev{X_{ij}}\,K\bigg(\frac{|X_{ij}|^{2}}{b^{2}}
    \bigg),
\end{align*}
where $X_{ij} = X_{j} - X_{i}$, $b$ is a \emph{bandwidth} and
$K(\cdot)$ is a univariate kernel supported on $[0,1]$. The
bandwidth $b$ should be selected so that the ball with the radius
$b$ and the center  at the point of estimation $X_{i}$ contains at
least  $d+1$ design points. For large value of $d$ this leads to a
bandwidth of order one and to a large estimation bias. The goal of
the structural adaptation is to diminish this bias using an
iterative procedure exploiting the available estimated structural
information.

In order to transform these general observations to a concrete
procedure, let us describe in the rest of this section how the
knowledge of the structure can help to improve the quality of the
estimation and how the structural information can be obtained when
the function or its estimator is given.

\subsection{Model estimation when an estimator of\/ $\cc{S}$ is available}
Let us start with the case of known $\cc{S}$.
The function $f$ has the same smoothness as $g$ in the
directions of the EDR subspace $\cc{S}$ spanned by
the vectors $\theta_{1},\ldots,\theta_{m^*}$, whereas it is constant
(and therefore, infinitely smooth) in all the orthogonal directions.
This suggests to apply an anisotropic bandwidth for
estimating the model function and its gradient. The corresponding
local-linear estimator can be defined by
\begin{align}
\label{hatfff}
\begin{pmatrix}
       \hat{f}(X_{i}) \\
       \widehat{\nablaf}(X_{i})
\end{pmatrix}
=\bigg\{
      \sum_{j=1}^{n}
        \onev{X_{ij}}
        \onev{X_{ij}}^{\T}
        w_{ij}^*
\bigg\}^{-1}
      \sum_{j=1}^{n}
      Y_{j} \onev{X_{ij}} \,
      w_{ij}^* \mbox{ } ,
\end{align}
with the weights $w_{ij}^* = K(|\Pi^* X_{ij}|^{2}/ h^{2})$, where
$h$ is some positive real number and $\Pi^*$ is the orthogonal
projector onto the EDR subspace $\cc{S}$.This choice of weights amounts
to using infinite bandwidth in the directions lying in the
orthogonal complement of the EDR subspace.

If only an estimator $\hat{\cc{S}}$ of $\cc{S}$ is available, the
orthogonal projector $\widehat\Pi$ onto the subspace $\hat{\cc{S}}$
may replace $\Pi^*$ in the expression (\ref{hatfff}). This rule of
defining the local neighborhoods is too stringent, since it
definitely discards the directions belonging to
$\hat{\cc{S}}^\perp$. Being not sure that our information about the
structure is exact, it is preferable to define the neighborhoods in
a softer way. This is done by setting
$w_{ij}=K(X_{ij}^\T(I+\rho^{-2}\widehat\Pi)X_{ij}/ h^{2})$ and by
redefining
\begin{align}
\label{hatfff1}
\begin{pmatrix}
       \hat{f}(X_{i}) \\
       \widehat{\nablaf}(X_{i})
\end{pmatrix}
=\bigg\{
      \sum_{j=1}^{n}
        \onev{X_{ij}}
        \onev{X_{ij}}^{\T}
        w_{ij}
\bigg\}^{-1}
      \sum_{j=1}^{n}
      Y_{j} \onev{X_{ij}} \,
      w_{ij} \mbox{ } .
\end{align}
Here, $\rho$ is a real number from the interval $[0,1]$ measuring the
importance attributed to the estimator $\widehat\Pi$. If we are very
confident
in our estimator $\widehat\Pi$, we should choose $\rho$ close to zero.

\subsection{Recovering the EDR subspace from an estimator of $\nablaf$}
Suppose first that the values of the function $\nablaf$ at the
points $X_{i}$ are known.Then $\cc{S}$ is the linear subspace of
$\RR^d$ spanned by the vectors $\nablaf(X_i)$, $i=1,\ldots,n$. For
classifying the directions of $\RR^d$ according to the variability
of $f$ in each direction and, as a by-product identifying $\cc{S}$,
the principal component analysis (PCA) can be used.

The PCA method is based on the orthogonal decomposition of the
matrix $\M = n^{-1} \sum_{i=1}^{n} \nablaf(X_{i})
\nablaf(X_{i})^{\T}$:  $\M = O \Lambda O^{T}$ with an
orthogonal matrix $O$ and a diagonal
matrix $\Lambda$ with diagonal entries $\lambda_{1} \ge \lambda_{2}
\ge \ldots \ge \lambda_{d}$. Clearly, for the multi-index model with
$m^*$-indices, only the first $m^*$ eigenvalues of $\M$ are
positive. The first $m^*$ eigenvectors of $\M$ (or, equivalently,
the first $m^*$ columns of the matrix $O$) define an orthonormal
basis in the EDR subspace .

Let $L$ be a positive integer. In \cite{HJPS}, a ``truncated''
matrix $\M_{L}$ is considered, which coincides with $\M$
if $L$ equals $n$. Let $\{ \psi_{\ell},
\ell=1,\ldots ,L \}$ be a system of functions on $\RR^{d}$ satisfying the
conditions
$
n^{-1} \sum_{i=1}^{n} \psi_{\ell}(X_{i}) \psi_{\ell'}(X_{i}) = \delta_{\ell,
\ell'}$ for every $\ell,\ell'\in\{1,\ldots,L\}$, with $\delta_{\ell,\ell'}$
being the Kronecker symbol.
Define
\begin{align}
\label{betal}
\beta_{\ell} = n^{-1} \sum_{i=1}^{n} \nablaf(X_{i}) \psi_{\ell}(X_{i})
\end{align}
and$\M_{L} = \sum_{\ell=1}^{L} \beta_{\ell} \beta_{\ell}^{\T}$.
By the Bessel inequality, it holds
$\M_{L} \le \M$. Moreover, since $\M \M_{L} = \M_{L} \M$,
any eigenvector of $\M$ is an eigenvector of $\M_{L}$.
Finally, by the Parseval equality, $\M_{L} = \M$ if
$L = n$. Note that the estimation of $\M$ has  been treated in \cite{Sam}.

The reason of considering the matrix $\M_L$ instead of $\M$ is that
$\M_L$ can be estimated much better than $\M$. In fact, estimators of $\M$
have poor performance for samples of moderate size because of the sparsity
of high dimensional data, ill-posedness of the gradient estimation and
the non-linear dependence of $\M$ on $\nabla f$. On the other hand,
estimation
of $\M_L$ reduces to the estimation of $L$ linear functionals of $\nabla f$
and may be done with a better accuracy. The obvious limitation of this
approach is that it recovers the EDR subspace entirely only if the rank
of $\M_L$ coincides with the rank of $\M$, which is equal to $m^*$.
To enhance our chances of seeing the condition $\rank(\M_L)=m^*$
fulfilled, we have to choose
$L$ sufficiently large. In practice, $L$ is chosen of the same order as $n$.

In the case when only an estimator of $\nabla f$ is available,
the above described method of recovering the EDR directions from
an estimator of $\M_L$ have a risk of order $\sqrt{L/n}$
(see \cite[Theorem 5.1]{HJPS}). This fact advocates against using
very large values of $L$. We desire nevertheless to use many
linear combinations in order to increase our chances of
capturing the whole EDR subspace. To this end, we modify the method of
extracting the structural information from the estimators $\hat\beta_\ell$
of vectors $\beta_\ell$.

Let $m\geq m^*$ be an integer.
Observe that the estimator $\widetilde\Pi_m$ of the projector $\Pi^*$
based on the PCA solves the following quadratic
optimization problem:
\begin{align}\label{PCA}
\text{minimize } \sum_{\ell}  \hat\beta_{\ell}^{\T}(I-\Pi) \hat\beta_{\ell}
\qquad
\text{subject to }
\qquad
\Pi^{2} = \Pi, \quad \tr \Pi \le m,
\end{align}
where the minimization is carried over the set of all symmetric $(d\times
d)$-matrices. The value $ m^{*}$ can be estimated by looking how many
eigenvalues of $\widetilde\Pi_m$ are significant.
Let $\cc{A}_{m}$ be the set of $ (d \times d) $-matrices defined as
follows:
\begin{align*}
\cc{A}_{m}
=
\{ \Pi: \Pi= \Pi^{\T} , \, 0 \preceq \Pi \preceq I, \,  \tr \Pi \le m \} .
\end{align*}
From now on, for two symmetric matrices $A$ and $B$, $A\preceq B$ means
that $B-A$ is semidefinite positive.
Define $\widehat{\Pi}_{m}$ as a minimizer of the  maximum
of the $\hat\beta_{\ell}^{\T}(I-\Pi) \hat\beta_{\ell}$'s instead of their
sum:
\begin{align}
\label{maxopt}
\widehat{\Pi}_{m} \in
\arginf_{\Pi \in \cc{A}_{m}} \max_{\ell} \hat\beta_{\ell}^{\T} (I-\Pi)
\hat\beta_{\ell}.
\end{align}
This is a convex optimization problem that can be effectively solved even
for
a large $d$ although a closed form solution is not known.
Moreover, as we will show below, the incorporation of (\ref{maxopt})
in the structural adaptation yields an algorithm having good
theoretical and empirical performance.

\section{Theoretical features of SAMM}
\label{Theor}

Throughout this section the true dimension $m^*$ of the EDR subspace is
assumed to be known. Thus, we are given $n$ observations
$(Y_1,X_1),\ldots,(Y_n,X_n)$  from the model
\begin{equation*}
Y_i=f(X_i)+\varepsilon_i=g(\theta_1^{\T}X_i,\ldots,\theta_{m^*}^{\T}X_i)
+\varepsilon_i,
\end{equation*}
where $\varepsilon_1,\ldots,\varepsilon_n$ are independent centered
random variables. The vectors $\theta_j$ are assumed to form an
orthonormal basis of the EDR subspace entailing thus the
representation $\Pi^*=\sum_{j=1}^{m^*}\theta_j\theta_j^{\T}$. In
what follows, we mainly consider deterministic design. Nevertheless,
the results hold in the case of random design as well, provided that
the errors are independent of $X_1,\ldots,X_n$. Henceforth, without
loss of generality we assume that $|X_i|<1$ for any $i=1,\ldots,n$,
where $|v|$ denotes the Euclidian norm of the vector $v$.

\subsection{Description of the algorithm}
The structure-adaptive algorithm with maximum minimization
consists of following steps.
\begin{itemize}
\item[a)] Specify positive real numbers $a_\rho$, $a_h$, $\rho_1$ and $h_1$.
Choose an integer $L$ and select a set $\{\psi_\ell,\,\ell\le L\}$
of vectors from $\RR^n$ verifying $|\psi_\ell|^2=n$. Set $k=1$.
\item[b)] Initialize the parameters $h=h_1$, $\rho=\rho_1$ and
$\widehat\Pi_0=0$.
\item[c)] Define the estimators $\widehat{\nablaf}(X_i)$ for $i=1,\ldots,n$
by formula (\ref{hatfff1}) with the current values of $h,\rho$ and
$\widehat\Pi$. Set
\begin{equation}\label{hbetal}
\hat\beta_\ell=\frac1n\sum_{i=1}^n \widehat{\nablaf}(X_i)\psi_{\ell,i},
\qquad \ell=1,\ldots,L,
\end{equation}
where $\psi_{\ell,i}$ is the $i$th coordinate of $\psi_\ell$.
\item[d)] Define the new value $\widehat\Pi_{k}$ as the solution to
(\ref{maxopt}).
\item[e)] Set $\rho_{k+1}=a_\rho\cdot\rho_k$, $h_{k+1}=a_h\cdot h_k$
and increase $k$ by one.
\item[f)] Stop if $\rho<\rho_{\rm{min}}$ or $h>h_{\rm{max}}$, otherwise
continue with the step c).
\end{itemize}
Let $k(n)$ be the total number of iterations. The matrix
$\widehat\Pi_{k(n)}$ is the desired estimator of the projector $\Pi^*$.
We denote by $\widehat\Pi_n$ the orthogonal projection onto the space spanned
by the eigenvectors of $\widehat\Pi_{k(n)}$ corresponding to the $m^*$
largest eigenvalues. The estimator of the EDR subspace is then
the image of $\widehat\Pi_n$. Equivalently, $\widehat\Pi_n$ is the estimator
of the projector onto $\cc{S}$.

The described algorithm requires the specification of the parameters
$\rho_1$, $h_1$, $a_\rho$ and $a_h$, as well as the choice of the set
of vectors $\{\psi_\ell\}$. In what follows we use the values
\begin{align*}
\begin{matrix}
\rho_1=1, &\rho_{\rm{min}}=n^{-1/(3\vee m^*)}, &a_\rho=e^{-1/{2(3\vee
m^*)}},\\
h_1=C_0n^{-1/(4\vee d)}, &h_{\rm{max}}= 2\sqrt{d}, &a_h=e^{1/2(4\vee d)}.
\end{matrix}
\end{align*}
This choice of input parameters is up to some minor modifications
the same as in \cite{HJS}, \cite{HJPS} and \cite{SSV}, and is based
on the trade-off between the bias and the variance of estimation. It
also takes into account the fact that the local neighborhoods used
in (\ref{hatfff}) should contain enough design points to entail the
consistency of the estimator. The choice of $L$ and that of vectors
$\psi_\ell$ will be discussed in Section~\ref{Sim}.

\subsection{Assumptions} Prior to stating rigorous theoretical results
we need to introduce a set of assumptions.
From now on, we use the notation $I$  for the identity matrix of dimension
$d$,
$\|A\|^2$ for the largest eigenvalue of $A^\T\cdot A$ and
$\|A\|_2^2$ for the sum of squares of all elements of the matrix
$A$.

\begin{description}
\item[(A1)] There exists a positive real $C_g$ such that $|\nabla g(x)|\le
C_g$
and $|g(x)-g(x')-(x-x')^\T \nablaG(x)|\le C_g|x-x'|^2$ for every
$x,x'\in\RR^{m^*}$.
\end{description}

Unlike the smoothness assumption, the assumptions on the
identifiability of the model and the regularity of design are more
involved and specific for each algorithm. The formal statements read
as follows.

\begin{description}
\item[(A2)] Let the vectors $\beta_\ell\in\RR^d$ be defined by (\ref{betal})
and let $\B =\bigl\{\betam = \sum_{\ell=1}^Lc_{\ell}\beta_\ell :
\sum_{\ell=1}^L |c_{\ell}|\le 1 \bigr\}$. There exist vectors
$\betam_{1},\ldots, \betam_{m^{*}} \in \B $
and constants $\mu_{1},\ldots ,\mu_{m^{*}}$ such that
\begin{align}
\label{Pi*}
\Pi^{*} \preceq \sum_{k=1}^{m^{*}} \mu_{k} \betam_{k} \betam_{k}^{\T} .
\end{align}
We denote $\mu^*=\mu_1+\ldots+\mu_{k}$.
\end{description}

\begin{remark} Assumption (A2) implies that the subspace $\cc{S}=Im(\Pi^*)$
is the smallest DR subspace, therefore it is the EDR
subspace. Indeed, for any DR subspace $\cc{S}'$, the gradient
$\nablaf(X_i)$ belongs to $\cc{S}'$ for every $i$. Therefore $\beta_\ell
\in\cc{S}'$ for every $\ell\leq  L$ and $\B\subset\cc{S}'$. Thus, for every
$\beta^\circ$ from the orthogonal complement ${\cc{S}'}^\perp$, it holds
$|\Pi^*\beta^\circ|^2\leq \sum_k \mu_k |\betam_k^\T \beta^\circ|^2=0$.
Therefore ${\cc{S}'}^\perp\subset\cc{S}^\perp$ implying thus the inclusion
$\cc{S}\subset\cc{S'}$.
\end{remark}

\begin{lemma}\label{lempsi} If the family $\{\psi_\ell\}$ spans $\RR^n$, then assumption (A2)
is always satisfied with some $\mu_k$ (that may depend on $n$).
\end{lemma}

\begin{Proof} Set $\Psi=(\psi_1,\ldots,\psi_L)\in\RR^{n\times L}$,
$\boldsymbol{\nablaf}=(\nablaf(X_1),\ldots,\nablaf(X_n))\in\RR^{d\times n}$ and
write the $d\times L$ matrix $B=(\beta_1,\ldots,\beta_L)$ in the form
$\boldsymbol{\nablaf}\cdot\Psi$.  Recall that if $M_1,M_2$ are two matrices such
that $M_1\cdot M_2$ is well defined and the rank of $M_2$ coincides with the
number of lines in $M_2$, then $\rank(M_1\cdot M_2)=\rank(M_1)$.  This implies
that $\rank(B)=m^*$ provided that $\rank(\Psi)=n$, which amounts to
$\Span(\{\psi_\ell\})=\RR^n$.

Let now $\tilde\beta_1,\ldots,\tilde\beta_{m^*}$ be a linearly independent
subfamily of $\{\beta_\ell,\ell\le L\}$.
Then the $m^*$th largest eigenvalue $\lambda_{m^*}(\tM)$ of the matrix
$\tM =\sum_{k=1}^{m^*}\tilde\beta_k\tilde\beta_k^\T$ is strictly
positive. Moreover, if $v_1,\ldots,v_{m^*}$ are the eigenvectors of
$\tM$ corresponding to the eigenvalues
$\lambda_1(\tM)\ge \ldots\ge \lambda_{m ^*}(\tM)>0$, then
$$
\Pi^*=\sum_{k=1}^{m^*} v_kv_k^\T\preceq
\frac1{\lambda_{m^*}}\sum_{k=1}^{m^*} \lambda_k v_kv_k^\T=
\lambda_{m^*}^{-1}\tM=
\lambda_{m^*}^{-1}\sum_{k=1}^{m^*}\tilde\beta_k\tilde\beta_k^\T.
$$
Hence, (\ref{Pi*}) is fulfilled with $\mu_k=1/\lambda_{m^*}(\tM)$
for every $k=1,\ldots,m^*$.
\end{Proof}

These arguments show that (A2) is a fairly weak identifiability
assumption. In fact, since we always choose $\{\psi_\ell\}$
so that $\Span(\{\psi_\ell\})=\RR^n$, (A2) amounts to requiring that
the value $\mu^*$ remains bounded when $n$ increases.

Let us proceed with the assumption on the design regularity.
Define $P_k^*=(I+\rho_k^{-2}\Pi^*)^{1/2}$,
$Z_{ij}^{(k)}=(h_kP_k^*)^{-1}X_{ij}$  and for any $d\times d$ matrix $U$
set $w_{ij}^{(k)}(U)=K\big((Z_{ij}^{(k)})^\T U Z_{ij}^{(k)}\big)$,
$\bar w_{ij}^{(k)}(U)=K'\big((Z_{ij}^{(k)})^\T U Z_{ij}^{(k)}\big)$,
$N_{i}^{(k)}(U)=\sum_j w_{ij}^{(k)}(U)$ and
$$
\V^{(k)}_{i}(U)=\sum_{j=1}^{n}
\onev{Z_{ij}^{(k)}}
\onev{Z_{ij}^{(k)}}^{\T}
\,w_{ij}^{(k)}(U).
$$

\begin{description}
\item[(A3)] For some positive constants $C_V,C_K,C_{K'},C_{w}$ and for some
$\alpha\in]0,1/2]$, the inequalities
\begin{align}
\|\V_i^{(k)}(U)^{-1}\|N_i^{(k)}(U)&\leq C_V,\qquad
i=1,\ldots,n,\label{assCV}\\
\sum_{i=1}^n w_{ij}^{(k)}(U)/N_i^{(k)}(U)&\leq C_K,\qquad
j=1,\ldots,n,\label{assCK}\\
\sum_{i=1}^n|\bar w^{(k)}_{ij}(U)|/N^{(k)}_i(U)&\leq C_{K'},\qquad
j=1,\ldots,n,\\
\sum_{j=1}^n |\bar w^{(k)}_{ij}(U)|/N_i^{(k)}(U)&\leq C_{w}\qquad
i=1,\ldots,n,\label{asCw}
\end{align}
hold for every $k\le k(n)$ and for every $d\times d$ matrix $U$ verifying
$\|U-I\|_2\le \alpha$.
\end{description}

\begin{description}
\item[(A4)] The errors $\{\varepsilon_i,i\le n\}$ are centered Gaussian with
variance $\sigma^2$.
\end{description}


\subsection{Main result}
We assume that the kernel $K$
used in (\ref{hatfff1}) is chosen to be continuous, positive and vanishing
outside the interval $[0,1]$. The vectors $\psi_\ell$ are assumed to verify
\begin{equation}
\label{asspsi}
\max_{\ell=1,\ldots,L}\max_{i=1,\ldots,n} |\psi_{\ell,i}|<\bar\psi,
\end{equation}
for some constant $\bar\psi$ independent of $n$.
In the sequel, we denote by $C,C_1,\ldots$
some constants depending only on $m^*,\mu^*,C_g,C_V,C_K,C_{K'},C_w$ and
$\bar\psi$.

\begin{theorem}\label{THM1}
Assume that assumptions \textbf{(A1)-(A4)} are fulfilled.
There exists a constant $C>0$ such that for any $z\in]0,2\sqrt{\log(nL)}]$
and for
sufficiently large values of $n$, it holds
$$
\Pb\bigg(\sqrt{\tr(I-\widehat\Pi_n)\Pi^*}>
Cn^{-\frac2{3\vee m^*}}t_n^2+\frac{2zc_0\sqrt{\mu^*}\sigma}{\sqrt{n(1-\zeta_n)}}
\bigg)\leq
Lze^{-\frac{z^2-1}2}+\frac{3k(n)-5}n,
$$
where $c_0=\bar\psi\sqrt{dC_KC_V}$, $t_n=O(\sqrt{\log(Ln)})$
and $\zeta_n=O(t_n\,n^{-\frac1{6\vee m^*}})$.
\end{theorem}

\begin{corollary}
Under the assumptions of Theorem~\ref{THM1}, for sufficiently large $n$,
it holds
\begin{align*}
&\Pb\bigg(\|\widehat\Pi_n-\Pi^*\|_2 >
Cn^{-\frac2{3\vee m^*}}t_n^2+\frac{2\sqrt{2\mu^*}zc_0\sigma}{\sqrt{n(1-\zeta_n)}}
\bigg)\leq
Lze^{-\frac{z^2-1}2}+\frac{3k(n)-5}n\\
&\esp(\|\widehat\Pi_n-\Pi^*\|_2 )\le C\bigg(n^{-2/(3\vee m^*)}t_n^2+
\frac{\sqrt{\log nL}}{\sqrt{n}}\bigg)+\frac{\sqrt{2m^*}(3k(n)-5)}n.
\end{align*}
\end{corollary}
\begin{Proof}
Easy algebra yields
\begin{align*}
\|\widehat\Pi_n-\Pi^*\|_2^2&=\tr(\widehat\Pi_n-\Pi^*)^2=\tr
\widehat\Pi_n^2-2\tr\widehat\Pi_n\Pi^*+\tr\Pi^*\\
&\leq \tr\widehat\Pi_n+m^*-2\tr\widehat\Pi_n\Pi^*\leq
2m^*-2\tr\widehat\Pi_n\Pi^*.
\end{align*}
The equality $\tr\Pi^*=m^*$ and the linearity of the trace operator
complete the proof of the first inequality. The second inequality
can be derived from the first one by standard arguments in view of
the inequality $\|\widehat\Pi_n-\Pi^*\|_2^2\le 2m^*$.
\end{Proof}

These results assess that for $m^*\leq 4$, the estimator of $\cc{S}$
provided by the SAMM procedure is $\sqrt{n}$-consistent up to a
logarithmic factor. This rate of convergence is known to be optimal
for a broad class of semiparametric problems, see \cite{BKRW} for a
detailed account on the subject.

\begin{remark}
The inspection of the proof of Theorem~\ref{THM1} shows that the
factor $t_n^2$ multiplying the ``bias'' term $n^{-2/(3\vee m^*)}$
disappears when $m^*>3$.
\end{remark}

\begin{remark}\label{Noise}
The same rate of convergence remains valid in the case when the errors
are not necessarily identically distributed Gaussian random variables,
but have (uniformly in $n$) a bounded exponential moment. This can be proved
along the lines of Proposition~\ref{prop1}, see Appendix.
\end{remark}

\begin{remark}
Note that in (A3) we implicitly assumed that the matrices $\V_i^{(k)}$
are invertible, which may be true only if any neighborhood
$E^{(k)}(X_i) = \{x : |(I+\rho_k^{-2}\Pi^*)^{-1/2}(X_i-x)|\le h_k\}$
contains
at least $d$ design points different from $X_i$.
The parameters $h_1$, $\rho_1$, $a_\rho$ and $a_h$ are chosen so that the
volume of ellipsoids $E^{(k)}(X_i)$ is a non-decreasing function of $k$
and $Vol(E^{(1)}(X_i))=C_0/n$. Therefore, from theoretical point
of view, if the design is random with positive density on $[0,1]^d$,
it is easy to check that for a properly chosen constant $C_0$, assumption
(A3) is satisfied with a probability close to one. In applications, we
define $h_1$ as the smallest real such that $\min_{i=1,\ldots,n}\#
E^{(1)}(X_i)=d+1$ and add to $\V_i$ a small full-rank matrix to be sure
that the resulting matrix is invertible, see Section \ref{Sim}.
\end{remark}

\begin{remark}
In the case when $m=n^{1/d}$ is integer,
an example of deterministic design satisfying (A3) is
as follows. Choose $d$ functions $h_k:[0,1]\to [0,\infty[$ such that
$\inf_{[0,1]} h_k(x)>0$
and $\sup_{[0,1]} h_k(x) <\infty$. Define the design points $\{X_i\}$
by $\{\int_{i_1/m}^{1+i_1/m}h_1(x)\,dx,\ldots
\int_{i_d/m}^{1+i_d/m}h_d(x)\,dx\}$, where $i_1,\ldots,i_d$ range over
$\{0,\ldots,m-1\}$. This definition guarantees that the number of design
points lying in an ellipsoid $E$ is asymptotically of the same order as
$n\,Vol(E)$, as $n\to\infty$. This suffices for (A3). Of course, it is
unlikely to have such a design in practice, since even for small $m$ and
moderate $d$ it leads to an unrealistically large sample size.
\end{remark}

\section{Simulation results}
\label{Sim}

The aim of this section is to demonstrate on several examples how the
performance of the algorithm SAMM depends on
the sample size $n$, the dimension $d$ and the noise level
$\sigma$. We also show that our procedure can be successfully
applied in autoregressive models. Many unreported results show
that in most situations
the performance of SAMM is comparable to the performance of SA
approach based on PCA and to that of MAVE. A thorough comparison
of the numerical virtues of these methods being out of scope of this
paper, we simply show on some examples that SAMM may substantially
outperform MAVE in the case of large ``bias''.

The computer code of the procedure SAMM is distributed freely, it can be
downloaded from
\textsl{\small http://www.proba.jussieu.fr/pageperso/dalalyan/}.
It requires the MATLAB packages SDPT3 and Yalmip.We are grateful
to Professor Yingcun Xia for making the computer code of MAVE
available to us.

To obtain higher stability of the algorithm, we preliminarily
standardize the response $Y$ and the predictors $X^{(j)}$. More
precisely, we deal with $\tilde Y_i=Y_i/\sigma_Y$ and $\tilde
X=\diag(\Sigma_X)^{-1/2}X$, where $\sigma_Y^2$ is the empirical
variance of $Y$, $\Sigma_X$ is the empirical covariance matrix of
$X$ and $\diag(\Sigma_X)$ is the $d\times d$ matrix obtained from
$\Sigma_X$ by replacing the off-diagonal elements by zero. To
preserve consistency, we set
${\tilde\beta}_{\ell,k(n)}=\diag(\Sigma_X)^{-1/2}\hat\beta_{\ell,k(n)}$,
where $\hat\beta_{\ell,k(n)}$ is the last-step estimate of
$\beta_\ell$, and define $\widehat\Pi_{k(n)}$ as the solution to
(\ref{maxopt}) with $\hat\beta_\ell$ replaced by
${\tilde\beta_{\ell,k(n)}}$. Furthermore, we add the small full-rank
matrix $I_{d+1}/n$ to $\sum_{j=1}^{n}\onev{X_{ij}}
\onev{X_{ij}}^{\T}w_{ij}$ in (\ref{hatfff1}).

In all examples presented below the number of replications is
$N=250$. The mean loss $\overline{{\text{er}}_N}=\frac1N\sum_j {\text{er}}_j$ and the
standard deviation $\sqrt{\frac1N\sum_j
({\text{er}}_j-\overline{{\text{er}}_N})^2}$ are reported, where
${\text{er}}_j=\|\widehat\Pi^{(j)}-\Pi^*\|$ with $\widehat\Pi^{(j)}$ being
the estimator of $\Pi^*$ for $j$th replication.

\subsection{Choice of $\{\psi_{\ell},\ell\le L\}$}
The set $\{\psi_\ell\}$ plays an essential role in the algorithm.
The optimal choice of this set is an important issue that needs
further investigation. We content ourselves with giving one particular
choice which agrees with theory and leads to nice empirical results.

Let $\mathfrak S_j$, $j\le d$, be the permutation of the set $\{1,\ldots,n\}$
satisfying $X_{\mathfrak S_j(1)}^{(j)}\le \ldots \le
X_{\mathfrak S_j(n)}^{(j)}$. Let $\mathfrak S_j^{-1}$ be the inverse of $\mathfrak S_j$, i.e.\
$\mathfrak S_j(\mathfrak S_j^{-1}(k))=k$ for every $k=1,\ldots,n$. Define $\{\psi_\ell\}$ as the
set of vectors
$$
\Bigg\{
\begin{matrix}
\Big(\cos\big( \frac{2\pi (k-1)\mathfrak S_j^{-1}(1)}{n}\big),\ldots,
\cos\big(\frac{2\pi (k-1) \mathfrak S_j^{-1}(n)}{n}\big)\Big)^\T\\[2pt]
\Big(\sin\big(\frac{2\pi k\mathfrak S_j^{-1}(1)}{n}\big),\ldots,
\sin\big(\frac{2\pi k\mathfrak S_j^{-1}(n)}{n}\big)\Big)^\T
\end{matrix},
k\le [n/2],\ j\le d\Bigg\}
$$
normalized to satisfy $\sum_{i=1}^n \psi_{\ell,i}^2=n$ for every $\ell$.
It is easily seen that these vectors satisfy conditions (\ref{asspsi}) and
$\Span(\{\psi_\ell\})=\RR^n$, so the conclusion of
Lemma~\ref{lempsi} holds. Above, $[n/2]$ is the integer part of $n/2$ and $k$ and $j$ are positive integers.

\subsection*{Example 1 (Single-index)}We set $d=5$ and
$f(x)=g(\theta^\T x)$ with
\begin{align*}
g(t)=4|t|^{1/2}\sin^2(\pi t),\quad\text{and}\quad
\theta=(1/\sqrt{5},2/\sqrt5, 0,0,0)^\T\in\RR^5.
\end{align*}
We run SAMM and MAVE procedures on the data generated by the model
$$
Y_i=f(X_{i})+0.5\cdot\varepsilon_i,
$$
where the design $X$ is such that the coordinates $(X_i^{(j)}, j\le
5,i\le n)$ are i.i.d.\ uniform on $[-1,1]$, and the errors
$\varepsilon_i$ are i.i.d.\ standard Gaussian independent of the
design.

\begin{figure} 
$
\begin{matrix}
\epsfxsize=120pt \epsfysize=110pt \epsffile{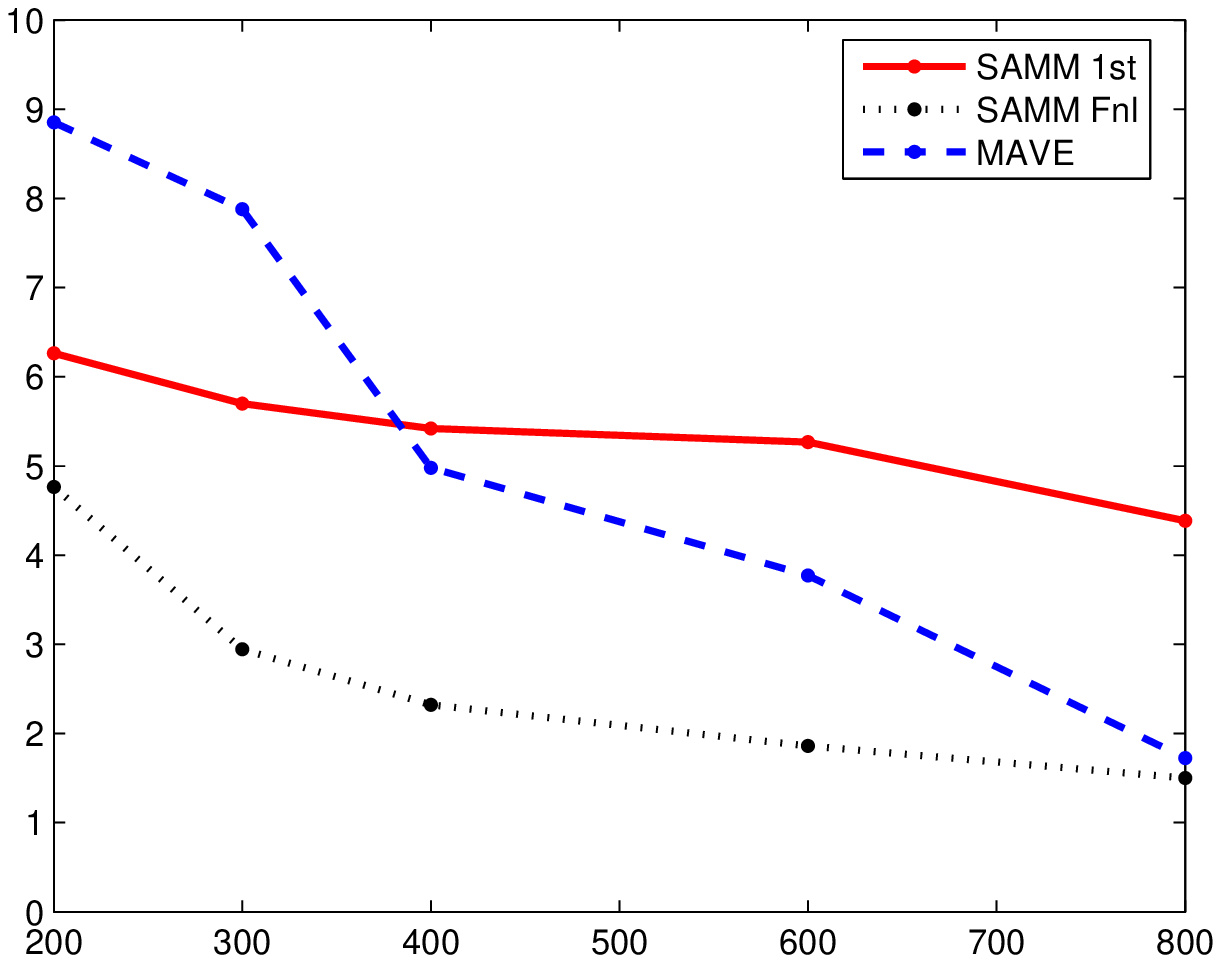}
&\text{\hglue-10pt\epsfxsize=120pt
\epsfysize=110pt\epsffile{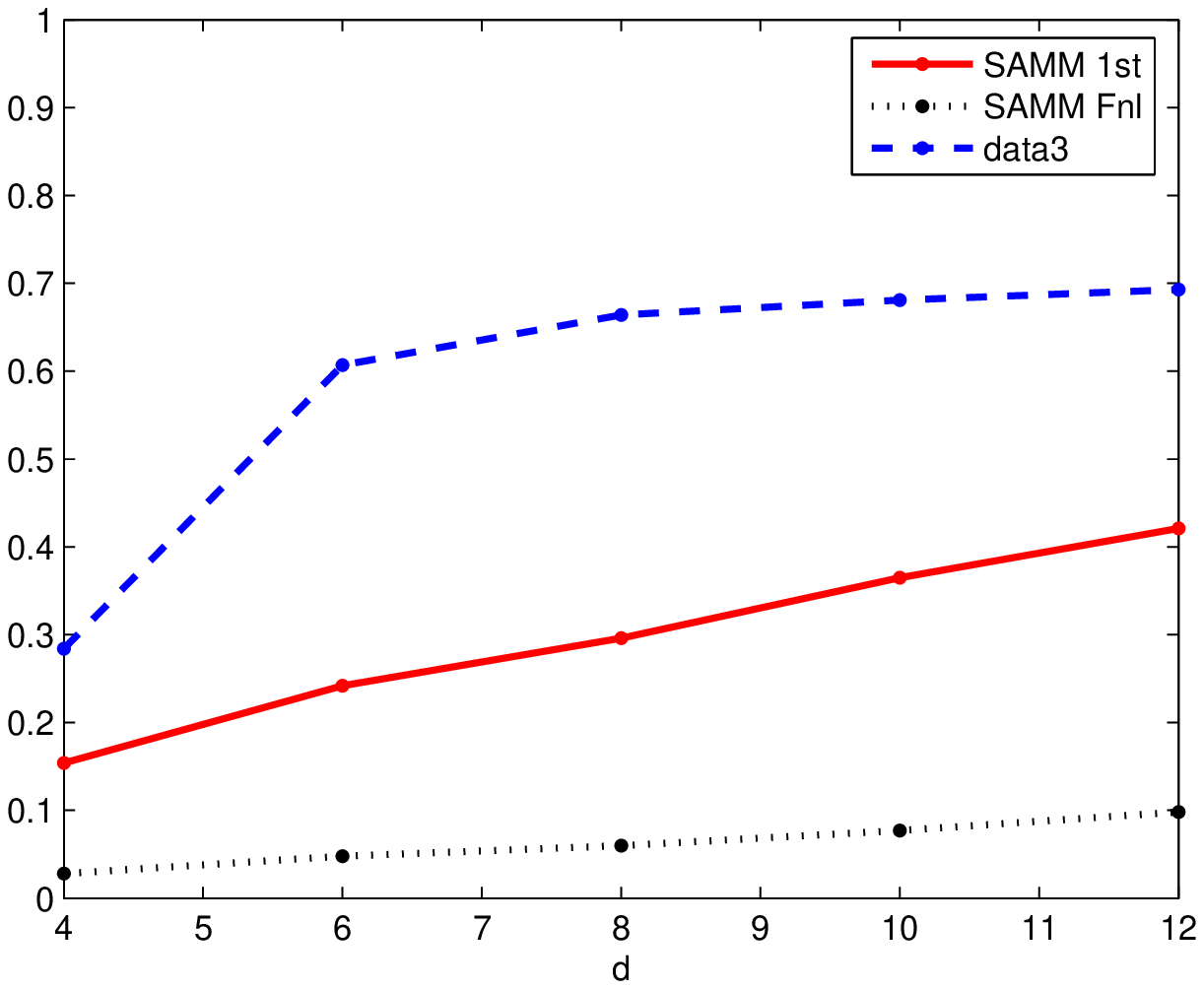}}
&\text{\hglue-10pt\epsfxsize=120pt \epsfysize=110pt\epsffile{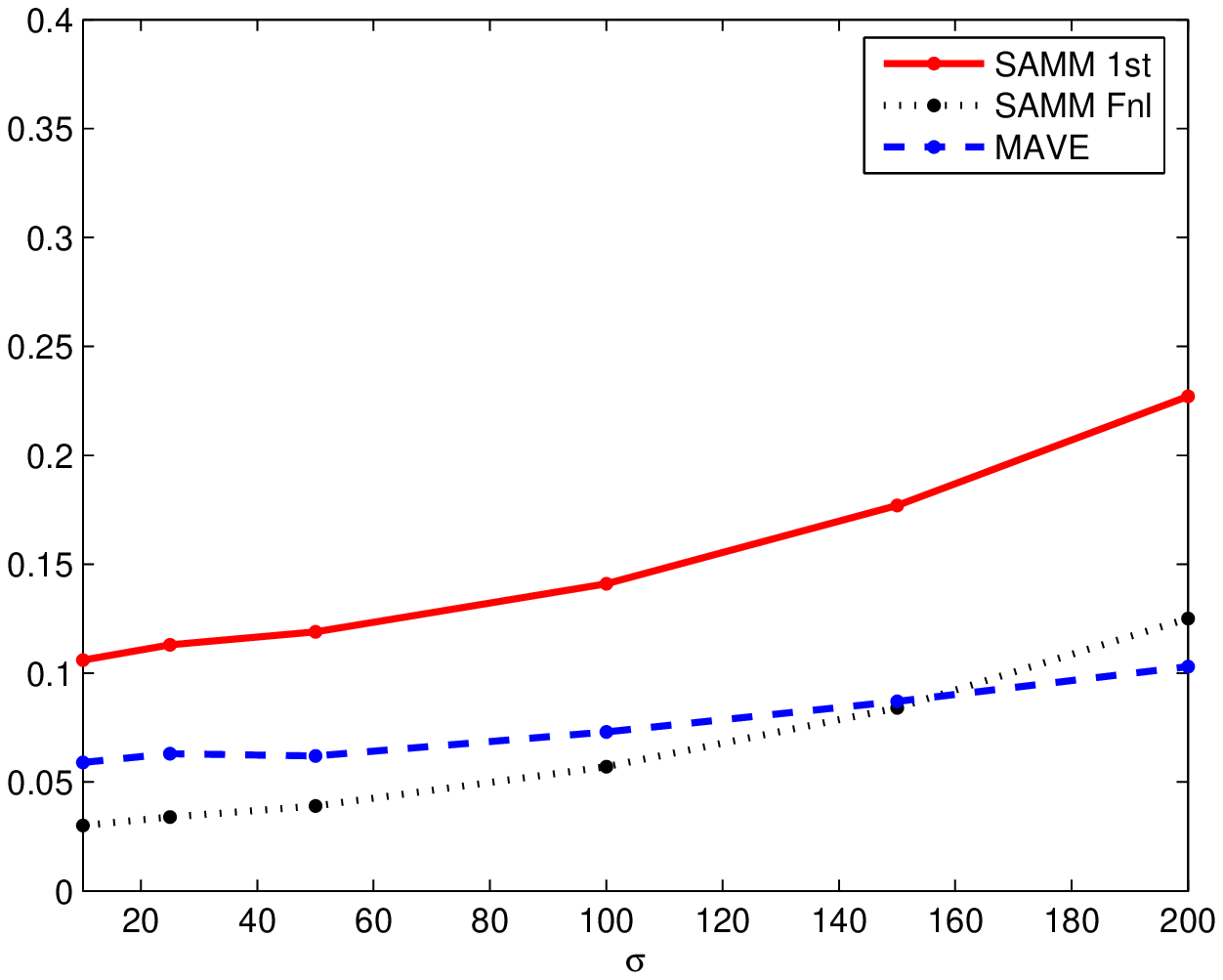}}\\
\text{\footnotesize\textbf{(a)}} & \text{\footnotesize\textbf{(b)}}&
\text{\footnotesize\textbf{(c)}}
\end{matrix}
$ \caption[]{\footnotesize (a) Average loss multiplied by $\sqrt{n}$
versus $n$ for the first step (full line) and the final (dotted
line) estimators provided by SAMM and for the estimator based on
MAVE (broken line) in Example 1, (b) (resp. (c)) Average loss versus
$d$ (resp. $\sigma$) for the first step (full line) and the final
(dotted line) estimators provided by SAMM and for the estimator
based on MAVE (broken line) in Example 2 (resp. Example 3).}
\label{fig1}
\end{figure}

Table \ref{table1} contains the average loss for different values of
the sample size $n$ for the first step estimator by SAMM, the final
estimator provided by SAMM and the estimator based on MAVE. We plot
in Figure~\ref{fig1} (a) the average loss normalized by the square
rood of the sample size $n$ versus $n$. It is clearly seen that the
iterative procedure improves considerably the quality of estimation
and that the final estimator provided by SAMM is
$\sqrt{n}$-consistent. In this example, MAVE method often fails to
recover the EDR subspace. However, the number of failures decreases
very rapidly with increasing $n$. This is the reason why the curve
corresponding to MAVE in Figure~\ref{fig1} (a) decreases with a
strong slope.

\begin{table}
\caption{\small\it Average loss $\|\widehat\Pi-\Pi^*\|$ of the
estimators obtained by SAMM and MAVE procedures in Example 1. The
standard deviation is given in parentheses.} \label{table1}
\medskip
\begin{tabular}{crrrrr}
\hline  $\boldsymbol{n}$ & 200 & 300 & 400 & 600 & 800\\
\hline  {\bf SAMM}, $\bf 1st$ & 0.443 & 0.329 & 0.271 & 0.215 & 0.155\\
& (.211) & (.120)& (.115) & (.095) &(.079) \\
 {\bf SAMM}, $\bf Fnl$ & 0.337 & 0.170 & 0.116 & 0.076 & 0.053\\
& (.273) & (.147)& (.104) & (.054) &(.031) \\
\bf  MAVE \phantom{$\bf 1st$} & 0.626 & 0.455 & 0.249 & 0.154 & 0.061\\
& (.363) & (.408)& (.342) & (.290) &(.161) \\
 \hline
\end{tabular}
\end{table}

\subsection*{Example 2} For $d\ge 2$ we set
$f(x)=g(\theta^\T x)$ with
$$
g(x)=(x_1-x_2^3)(x_1^3+x_2);
$$
and  $\theta_1=(1,0,\ldots,0)\in \RR^d$,
$\theta_2=(0,1,\ldots,0)\in\RR^d$. We run SAMM and MAVE procedures
on the data generated by the model
$$
Y_i=f(X_{i})+0.1\cdot\varepsilon_i,\qquad i=1,\ldots,300,
$$
where the design $X$ is such that the coordinates $(X_i^{(j)}, j\le
d,i\le n)$ are i.i.d.\ uniform on $[-40,40]$, and the errors
$\varepsilon_i$ are i.i.d.\ standard Gaussian independent of the
design. The results of simulations for different values of $d$ are
reported in Table~\ref{table2}.

As expected, we found that (cf.\ Figure~\ref{fig1}(b))  the quality
of SAMM deteriorated linearly in $d$ as $d$ increased. This agrees
with our theoretical results. It should be noted that in this case
MAVE fails to find the EDR space.

\begin{table}
\caption{\small\it Average loss $\|\widehat\Pi-\Pi^*\|$ of the
estimators obtained by SAMM and MAVE procedures`in Example 2. The
standard deviation is given in parentheses.} \label{table2}
\medskip
\begin{tabular}{crrrrr}
\hline  $\boldsymbol{d}$ & 4 & 6 & 8 & 10 & 12\\
\hline  \bf  SAMM $\bf 1st$ & 0.154 & 0.242 & 0.296 & 0.365 & 0.421\\
& (.063) & (.081)& (.071) & (.087) &(.095) \\
 {\bf SAMM}, $\bf Fnl$ & 0.028 & 0.048 & 0.060 & 0.077 & 0.098\\
& (.011) & (.020)& (.021) & (.026) &(.037) \\
\bf  MAVE \phantom{$\bf 1st$} & 0.284 & 0.607 & 0.664 & 0.681 & 0.693\\
& (.147) & (.073)& (.052) & (.054) &(.044) \\
 \hline
\end{tabular}
\end{table}

\subsection*{Example 3}
For $d=5$ we set $f(x)=g(\theta^\T x)$ with
$$
g(x)=(1+x_1)(1+x_2)(1+x_3)
$$
and  $\theta_1=(1,0,0,0,0)$, $\theta_2=(0,1,0,0,0)$,
$\theta_3=(0,0,1,0,0)$. We run SAMM and MAVE procedures on the data
generated by the model
$$
Y_i=f(X_{i})+\sigma\cdot\varepsilon_i,\qquad i=1,\ldots,250,
$$
where the design $X$ is such that the coordinates $(X_i^{(j)}, j\le
d,i\le n)$ are i.i.d.\ uniform on $[0,20]$, and the errors
$\varepsilon_i$ are i.i.d.\ standard Gaussian independent of the
design.

Figure \ref{fig1}(c) shows that the qualities of both SAMM and MAVE
deteriorate linearly in $\sigma$, when $\sigma$ increases. These
results also demonstrate that, thanks to an efficient bias
reduction, the SAMM procedure outperforms MAVE when stochastic error
is small, whereas MAVE works better than SAMM in the case of
dominating stochastic error (that is when $\sigma$ is large).

\begin{table}
\caption{\small\it Average loss $\|\widehat\Pi-\Pi^*\|$ of the
estimators obtained by SAMM and MAVE procedures in Example 3. The
standard deviation is given in parentheses.} \label{table3}
\medskip
\begin{tabular}{crrrrrr}
\hline  $\boldsymbol{\sigma}$ & 200 & 150 & 100 & 50 & 25 & 10\\
\hline  \bf  SAMM $\bf 1st$ & 0.227 & 0.177 & 0.141 & 0.119 & 0.113 & 106\\
& (.092) & (.075)& (.055) & (.051) &(.048) & (.043) \\
 {\bf SAMM}, $\bf Fnl$ & 0.125 & 0.084 & 0.057 & 0.039 & 0.034 & 0.03 \\
& (.076) & (.037)& (.026) & (.019) &(.021) & (.018) \\
\bf  MAVE \phantom{$\bf 1st$} & 0.103 & 0.087 & 0.073 & 0.062 & 0.063 & 0.059\\
& (.041) & (.035)& (.027) & (.023) &(.024) & (.023) \\
 \hline
\end{tabular}
\end{table}

\subsection*{Example 4 (time series)}
Let now $T_1,\ldots,T_{n+6}$ be generated by the autoregressive
model
$$
T_{i+6}=f(T_{i+5},T_{i+4},T_{i+3},T_{i+2},T_{i+1},T_i)+0.2\cdot\varepsilon_i,\qquad
i=1,\ldots,n,
$$
with initial variables $T_1,\ldots,T_6$ being independent standard
normal independent of the innovations $\varepsilon_i$, which are
i.i.d.\ standard normal as well. Let now $f(x)=g(\theta^\T x)$ with
\begin{align*}
g(x)&=-1+0.6x_1-\cos(0.5\pi x_2)+e^{-x_3^2},\\
\theta_1&=(1,0,0,2,0,0)/\sqrt5,\\
\theta_2&=(0,0,1,0,0,2)/\sqrt5,\\
\theta_3&=(-2,2,-2,1,-1,1)/\sqrt{15}.
\end{align*}
We run SAMM and MAVE procedures on the data $(X_i,Y_i)$,
$i=1,\ldots,250$, where $Y_i=T_{i+6}$ and
$X_i=(T_i,\ldots,T_{i+5})^\T$. The results of simulations reported
in Table~\ref{table4} show that the qualities of SAMM and MAVE are
comparable, with SAMM being slightly more performant.

\begin{table}
\caption{\small\it Average loss $\|\widehat\Pi-\Pi^*\|$ of the
estimators obtained by SAMM and MAVE procedures in Example 4. The
standard deviation is given in parentheses.} \label{table4}
\medskip
\begin{tabular}{crrrr}
\hline  $\boldsymbol{n}$ & 300 & 400 & 500 & 600\\
\hline  {\bf SAMM}, $\bf 1st$ & 0.391 & 0.351 & 0.334 & 0.293 \\
& (.172) & (.161)& (.137) & (.132) \\
 {\bf SAMM}, $\bf Fnl$ & 0.220 & 0.186 & 0.174 & 0.146 \\
& (.119) & (.123)& (.102) & (.089) \\
\bf  MAVE \phantom{$\bf 1st$} & 0.268 & 0.231 & 0.209 & 0.182 \\
& (.209) & (.170)& (.159) & (.122) \\
 \hline
\end{tabular}
\end{table}

\section*{Appendix}
Since the proof of the main result is carried out in several steps,
we give a short road map for guiding the reader throughout the
proof. The main idea is to evaluate the accuracy of the first step
estimators of $\beta_\ell$ and, given the accuracy of the estimator
at the step $k$, evaluate the accuracy of the estimators at the step
$k+1$. This is done in Subsections~\ref{ss2} and \ref{ss3}. These
results are based on a maximal inequality proved in Subsection
\ref{ss5} and on some properties of the solution to (\ref{maxopt})
proved in Subsection \ref{ss6}. The proof of Theorem~\ref{THM1} is
presented in Subsection \ref{ss4}, while some technical lemmas are
postponed to Subsection~\ref{ss7}.

\subsection{The accuracy of the first-step estimator}\label{ss2}

Since at the first step no information about the EDR subspace is
available, we use the same bandwidth in all directions, that is the
local neighborhoods are balls (and not ellipsoids) of radius $h$.
Therefore the first step estimator $\hat\beta_{1,\ell}$ of the
vector $\beta^*_\ell$ is the same as the one used in \cite{HJPS}.

\begin{proposition}\label{P6}
Under assumptions (A1),(A3), (A4) and (\ref{asspsi}), for every $\ell\le L$,
$$
|\hat\beta_{1,\ell}-\beta_{\ell}|\leq h_1C_g\sqrt{2C_V}
+\frac{\xi_{1,\ell}}{h_1\sqrt n},
$$
where $\xi_{1,\ell}$ is a zero mean normal vector verifying
$\esp|\xi_{1,\ell}|^2\leq
2d\sigma^2 C_V C_K \bar\psi^2$.
\end{proposition}

\begin{Proof}
Since at the first iteration we take $S_1=I$, the inequality
$|S_1X_{ij}|\leq h_1$ implies that
$|\Pi^* X_{ij}|\leq |X_{ij}|\leq h_1$. Therefore the bias term
$|P_1^*(E\hat\beta_{1,\ell}-\beta_{\ell})|$ is bounded by
$h_1C_g\sqrt{C_V}$ (cf.\ the proof of
Proposition~\ref{P5}).

For the stochastic term, we set $\xi_{1,\ell}=h_1\sqrt
n(\hat\beta_{1,\ell}-
E\hat\beta_{1,\ell})$. By Lemma~\ref{lem4.1}, we have
$\esp|P^*_1\xi_{1,\ell}|^2\leq d\sigma^2 C_V C_K\bar\psi^2$.
The assertion of the proposition follows now from
$P_1^*=(I+\rho^{-2}_1\Pi^*)^{-1/2}\succeq I/\sqrt2$.
\end{Proof}

\begin{corollary}\label{cor2} If\/ $nL\ge 6$ and the assertions of
Proposition~\ref{P6} hold, then
$$
\Pb\bigg(\max_\ell|\hat\beta_{1,\ell}-\beta_{\ell}|\geq h_1C_g\sqrt{C_V}+
\frac{2\sqrt{2dC_VC_K\log(nL)}\,\sigma\bar\psi}{h_1\sqrt n}\bigg)\leq
\frac1n.
$$

\end{corollary}

\begin{remark}
In order that the kernel estimator of $\nabla f(x)$ be consistent,
the ball centered at $x$ with radius $h_1$ should contain at least
$d$ points from $\{X_i,i=1,\ldots,n\}$. If the design is regular,
this means that $h_1$ is at least of order $n^{-1/d}$. The
optimization of the risk of $\hat\beta_{1,\ell}$ with respect to
$h_1$ verifying $h_1\geq n^{-1/d}$ leads to the choice
$h_1=Const.n^{-1/(4\vee d)}$.
\end{remark}

\subsection{One step improvement}\label{ss3}
At the $ k $th step of iteration, we have at our disposal a symmetric
matrix
$ \Pi\in\mathcal M_{d\times d} $ belonging to the set
$$
\mathscr P_{\delta}(\Pi^*)=\Big\{\Pi\in\mathcal M_{d\times d}:\
\tr\Pi\leq m^*,\quad 0\preceq \Pi\preceq I,\quad \tr(I-\Pi)\Pi^{*}\leq
\delta^{2}\Big\}
$$
Thus the matrix $ \Pi $ is the $ k $th step approximation of the
projector $ \Pi^{*} $ onto the EDR subspace $ \cc{S}^{*} $. Using
this
approximation, we construct the new matrix $ \hat \Pi $ in the
following way: Set $ S_{\Pi,\rho}=(I+\rho^{-2}\Pi)^{1/2} $,
$ P_{\Pi,\rho}=S_{\Pi,\rho}^{-1} $ and define the estimator of the
regression function and its gradient at the design point $ X_{i} $ as
follows:
$$
\begin{pmatrix}
\hat f_{\Pi}(X_{i})\\[3pt]
\widehat{\nablaf}_{\Pi}(X_{i})
\end{pmatrix}
=
V_{i}(\Pi)^{-1}\sum_{j=1}^{n} Y_{j}
\onev{X_{ij}}
w_{ij}(\Pi),
$$
where $ w_{ij}(\Pi)=K\bigl(h^{-2}|S_{\Pi,\rho}X_{ij}|^{2}\bigr) $ and
$$
V_{i}(\Pi)=\sum_{j=1}^{n}
\onev{X_{ij}}
\onev{X_{ij}}^{\T}
\,w_{ij}(\Pi).
$$
To state the next result, we need some additional notation. Set
$Z_{ij}=(hP_\rho^*)^{-1}X_{ij}$,
$U=P_\rho^*S_{\Pi,\rho}^2 P_\rho^*$ and $U^*=I$, where
$P_{\rho}^{*}=P_{\Pi^{*},\rho}=(I-\Pi^{*})+\rho(1+\rho^{2})^{-1/2}\;\Pi^{*}$.
In this notation, we obtain
\begin{align*}
\begin{pmatrix}
h^{-1}\hat f_{\Pi}(X_{i})\\[3pt]
P_\rho^*\widehat{\nablaf}_{\Pi}(X_{i})
\end{pmatrix}
&= \begin{pmatrix}
    h^{-1} & 0\\
    0 & P_{\rho}^{*}
\end{pmatrix}
V_{i}(\Pi)^{-1}\sum_{j=1}^{n} Y_j \onev{X_{ij}}\,w_{ij}(\Pi)\\
&= \frac1h
\begin{pmatrix}
    1 & 0\\
    0 & hP_{\rho}^{*}
\end{pmatrix}
V_{i}(\Pi)^{-1}\begin{pmatrix}
    1 & 0\\
    0 & hP_{\rho}^{*}
\end{pmatrix}
\sum_{j=1}^{n} Y_j \onev{Z_{ij}}\,w_{ij}(\Pi)\\
&= h^{-1}
\V_{i}(U)^{-1}
\sum_{j=1}^{n} Y_j \onev{Z_{ij}}\,w_{ij}(U)
\end{align*}
where $ w_{ij}(U)=K\big(Z_{ij}^\T UZ_{ij}\big)$ and
$$
\V_{i}(U)=\sum_{j=1}^{n}
\onev{Z_{ij}}
\onev{Z_{ij}}^{\T}
\,w_{ij}(U).
$$
Set $N_i(U)=\sum_j w_{ij}(U)$ and
$\alpha=2\delta^2\rho^{-2}+2\delta\rho^{-1}$.

\begin{proposition}\label{P5}
If (A1)-(A4) are fulfilled
then there exist Gaussian vectors $\xi_1^*,\ldots,\xi_L^*\in\RR^d$ such that
$\esp[|\xi_\ell^*|^2]\leq c_0^2\sigma^2$ and
\begin{align*}
\Pb\bigg(\sup_{\Pi,\ell}
\Big|P_\rho^*(\hat\beta_{\ell,\Pi}-\beta_{\ell})-\frac{\xi_\ell^*}{h\sqrt
n}\Big|
\geq \sqrt C_V C_g(\rho+\delta)^2h+\frac{c_1\sigma\alpha
t_n}{h\sqrt{n}}\bigg)\leq \frac2n,
\end{align*}
where the $\sup$ is taken over $\Pi\in\mathscr P_{\delta}$,
$\ell=1,\ldots,L$ and we used the notation
$t_n=5+\sqrt{3\log(Ln)+\frac32d^2\log n}$,
$c_0=\bar\psi\sqrt{dC_KC_V}$ and
$c_1= 30\bar\psi(C_w^2C_V^4C_K^2+C_V^2C_{K'}^2)^{1/2}$.
\end{proposition}

\begin{Proof}
Let us start with evaluating the bias term $|P_\rho^*(\esp \hat\beta_{\ell,\Pi}
-
\beta_{\ell})|$.
According to the Cauchy-Schwarz inequality, it holds
\begin{align*}
\big|P_{\rho}^{*}\bigl(\esp \hat\beta_{\ell,\Pi} - \beta_{\ell}\bigr)\big|^2
&=
n^{-2}
\bigg|
    \sum_{i=1}^{n}
P_{\rho}^{*}\bigl(\esp[\widehat{\nablaf}_{\Pi}(X_{i})]-\nablaf(X_{i})\bigr)
    \psi_{\ell}(X_{i})
\bigg|^2
\nonumber\\
&\leq
\frac1{n^{2}}
    \sum_{i=1}^{n}
\big|P_{\rho}^{*}\bigl(\esp[\widehat{\nablaf}_{\Pi}(X_{i})]-
    \nablaf(X_{i})\bigr)\big|^2
    \sum_{i=1}^{n} \psi_{l}^{2}(X_{i})
\nonumber\\
&\leq
\max_{i=1,\ldots,n}\big|P_{\rho}^{*}\bigl(\esp[\widehat{\nablaf}_{\Pi}(X_{i})]-
\nablaf(X_{i})\bigr)\big|^2.
\end{align*}
Simple computations show that
\begin{align*}
\big|P_{\rho}^{*}\bigl(\esp[\widehat{\nablaf}_{\Pi}(X_{i})]&-
\nablaf(X_{i})\bigr)\big|\\
&\leq
\Bigg|\esp
\begin{pmatrix}
    h^{-1}\hat f_{\Pi}(X_{i})\\[3pt]
    P_{\rho}^{*}\widehat{\nablaf}_{\Pi}(X_{i})
\end{pmatrix}
-
\begin{pmatrix}
    h^{-1}f(X_{i})\\[3pt]
    P_{\rho}^{*}\nablaf(X_{i})
\end{pmatrix}\Bigg|\\
&=h^{-1}\Bigg|
\V_{i}^{-1}
\sum_{j=1}^{n} f(X_j) \onev{Z_{ij}}\,w_{ij}(U)-
\begin{pmatrix}
    h^{-1}f(X_{i})\\[3pt]
    P_{\rho}^{*}\nablaf(X_{i})
\end{pmatrix}\Bigg|\\
&=h^{-1}
\Big|\V_{i}^{-1}
\sum_{j=1}^{n} r_{ij} \onev{Z_{ij}} w_{ij}(U)\Big|:=b(X_i),
\end{align*}
where $ r_{ij} = f(X_{j}) - f(X_{i}) - X_{ij}^{\T}\nablaf(X_{i}) $.
Define
$ \lambda_{j}=h^{-1}r_{ij}\sqrt{w_{ij}(U)} $ and
$v_{j} = \V_{i}^{-1/2}\onev{Z_{ij}}\sqrt{w_{ij}(U)}$.
Then
$$ b(X_i) =
\biggl| \V_{i}^{-1/2}\sum_{j=1}^{n}\lambda_{j} v_{j} \biggr|
\le \bigl\| \V_{i}^{-1/2} \bigr\| \cdot |\lambda| \cdot
\biggl\| \sum_{j=1}^n v_{j} v_{j}^{\T} \biggr\|^{1/2}.
$$
The identity
$ \sum_{j} v_{j} v_{j}^{\T} = I_{d+1} $ implies
\begin{align*}
b(X_{i})^{2}
&\leq
\frac{1}{h^{2}} \;
\left\| \V_{i}^{-1/2}\right\|^{2} \cdot \sum_{j=1}^{n} r_{ij}^{2}
w_{ij}(U)
\nonumber\\
&\leq
h^{-2}\max_{j} r_{ij}^{2}\;
\left\| \V_{i}^{-1} \right\| \cdot\sum_{j=1}^{n} w_{ij}(U)
\nonumber\\
&\leq
C_V h^{-2}\max_{j} r_{ij}^{2} \, ,
\end{align*}
where the maximum of $ r_{ij} $ is taken over the indices $ j $
satisfying $ w_{ij}(U)\not=0 $. Since the weights $ w_{ij} $ are
defined via the kernel function $ K $ vanishing on the interval
$ [1,\infty[ $, we have$\max_{j} r_{ij}=\max\{r_{ij}: |S_{\Pi,\rho}X_{ij}|
\leq h\}$.
By Corollary~\ref{Cellipt} $ |S_{\Pi,\rho}X_{ij}|\leq h $ implies
$ |\Pi^{*}X_{ij}|\leq (\rho+\delta)h $. Let us denote by
$\Theta$ the $(d\times m^*)$ matrix having $\theta_k$ as $k$th column.
Then $\Pi^*=\Theta\Theta^{\T}$ and therefore
\begin{align*}
|r_{ij}|
&=
|f(X_{j})-f(X_{i})-X_{ij}^{\T}\nablaf(X_{i})|
\nonumber\\
&=
|g(\Theta^\T X_{j})-g(\Theta^\T X_{i})-(\Theta^\T X_{ij})^{\T}
\nablaG(\Theta^\T X_{i})|\\
&\le  C_{g} |\Theta^\T X_{ij}|^{2}\leq C_{g}(\rho+\delta)^{2} h^{2}.
\end{align*}
These estimates yield
$ |b(X_{i})|\leq \sqrt{C_V} \,C_{g}(\rho+\delta)^{2}h $, and consequently,
\begin{align}
\big|P_{\rho}^{*}\bigl(\esp \hat\beta_{\ell,\Pi} - \beta_{\ell}\bigr)\big|
\leq \max_i b(X_i)\leq
\sqrt{C_V} \, C_g (\rho+\delta)^{2} h.
\end{align}
Let us treat now the stochastic term
$P_{\rho}^{*}\bigl(\hat\beta_{\ell,\Pi}-\beta_{l}^*\bigr)$.
It can be bounded as follows
\begin{align*}
\big|
    P_{\rho}^{*}
    \bigl( \hat\beta_{\ell,\Pi}-\esp\hat\beta_{\ell,\Pi} \bigr)
\big|
\leq
\bigg|
    \sum_{j=1}^{n} c_{j,\ell}(U)\,\varepsilon_{j}
\bigg|,
\end{align*}
where
$$
c_{j,\ell}(U)=\frac{1}{hn}\sum_{i=1}^{n} \V_{i}^{-1}(U) \onev{Z_{ij}}
w_{ij}(U)\psi_{\ell}(X_{i}).
$$
Let us define $\xi_\ell^*=h\sqrt n\,
P_\rho^*(\hat\beta_{\ell,\Pi^*}-\esp[\hat\beta_{\ell,\Pi^*}])$.
In view of Lemma~\ref{lem4.1}, we have $\esp[|\xi_\ell^*|^2]\leq
nh^2\sigma^2\sum_j|c_{j,\ell}(U^*)|^2\leq c_0^2\sigma^2$.

One checks that for any $\ell=1,\ldots, L$ and for any
$\Pi$ such that $\tr(I-\Pi)\Pi^*\leq \delta^2$, it holds
$$
\Big|P_\rho^*(\hat\beta_{\ell,\Pi}-\esp[\hat\beta_{\ell,\Pi}])
-\frac{\xi_\ell^*}{h\sqrt n}\Big|\leq
\sup_{\|U-U^*\|_2\leq \alpha}
\bigg|
    \sum_{j=1}^{n} \big(c_{j,\ell}(U)-c_{j,\ell}(U^*)\big)\,\varepsilon_{j}
\bigg|.
$$
Set $a_{j,\ell}(U)=c_{j,\ell}(U)-c_{j,\ell}(U^*)$. Lemma~\ref{lem4.2}
implies that
Proposition~\ref{prop1} can be applied with
$\kappa_0=\frac{c_1\alpha}{h\sqrt{n}}$ and $\kappa_1=\frac{c_1}{h\sqrt n}$.
Setting $\epsilon=2\alpha/\sqrt n$ we get that the probability of the event
\begin{equation*}
\bigg\{\sup_{U,\ell}\bigg|\sum_{j=1}^{n}
\big(c_{j,\ell}(U)-c_{j,\ell}(U^*)\big)\,\varepsilon_{j}
\bigg|\geq \frac{c_1\sigma\alpha(5+\sqrt{3\log(Ln)+3d^2\log(\sqrt
n)})}{h\sqrt n} \bigg\}
\end{equation*}
is less than $2/n$. This completes the proof of the proposition.
\end{Proof}

\begin{corollary}\label{cor1}
If $nL\geq 6$ and the assumptions of Proposition~\ref{P5} are fulfilled,
then
\begin{align*}
\Pb\bigg(\sup_{\ell,\Pi}
\big|P_\rho^*(\hat\beta_{\ell,\Pi}-\beta_{\ell})\big|
\geq \sqrt C_V C_g(\rho+\delta)^2h+\frac{\sigma(zc_0+c_1\alpha t_n)}{h\sqrt
n}
\bigg)\leq Lze^{-\frac{z^2-1}2}.
\end{align*}
In particular, if $nL\ge 6$, the probability of the event
\begin{align*}
\bigg\{\sup_{\ell,\Pi}
\big|P_\rho^*(\hat\beta_{\ell,\Pi}-\beta_{\ell})\big|
\geq \sqrt C_V
C_g(\rho+\delta)^2h+\frac{\sigma(2c_0\sqrt{\log(Ln)}+c_1\alpha t_n)}{h\sqrt
n}
\bigg\}
\end{align*}
does not exceed $3/n$,
where $\sup$ is taken over all $\Pi\in\mathscr P_{\delta}(\Pi^*)$,
$\ell=1,\ldots,L$ and $c_0$, $c_1$, $t_n$
are defined in Proposition~\ref{P5} and in Theorem~\ref{THM1}.
\end{corollary}
\begin{Proof}
In view of Lemma 7 in \cite{HJS} and Lemma~\ref{lem4.1} ,
we have
$$
\Pb\big(\max_{\ell=1,\ldots,L}\big|\xi_\ell^*\big|\geq z c_0\sigma\big)\leq
\sum_{\ell=1}^L\Pb\big(\big|\xi_\ell^*\big|\geq z c_0\sigma\big)\leq L
ze^{-(z^2-1)/2}.
$$
The choice $z=\sqrt{4\log(nL)}$ leads to the desired inequality
provided that $nL\geq 6$.
\end{Proof}

\subsection{Proof of Theorem~\ref{THM1}}\label{ss4}

Recall that at the first step we use the following values of parameters:
$\widehat\Pi_0=\boldsymbol{0}$, $\rho_1=1$
and $h_1=n^{-1/(d\vee 4)}$. Let us denote
\begin{align*}
\gamma_1&= h_1C_g\sqrt C_V
+\frac{2\sqrt{2dC_VC_K\log(nL)}\,\sigma\bar\psi}{h_1\sqrt n},\quad
\delta_1= 2\gamma_1\sqrt{\mu^*},
\end{align*}
and introduce the event
$\Omega_1=\{\max_\ell|\hat\beta_{1,\ell}-\beta_{\ell}|\leq \gamma_1\}$.
According to Corollary~\ref{cor2} the probability
of the event $\Omega_1$ is at least $1-n^{-1}$. In view of
Proposition~\ref{P2}, we get
$\Pb(\tr(I-\widehat\Pi_1)\Pi^*\le \delta^2_1)\ge 1-n^{-1}$.

For any integer $k\in [2,k(n)]$ (where $k(n)$ is the total number of
iterations), we define
\begin{align*}
\rho_k&= a_\rho\rho_{k-1},\quad h_k=a_hh_{k-1},\quad
\alpha_k=\frac{2\delta_{k-1}}{\rho_k}\bigg(\frac{\delta_{k-1}}{\rho_k}+1\bigg),\\
\gamma_k&= \begin{cases}
\displaystyle C_g\sqrt
C_V\,(\rho_k+\delta_{k-1})^2h_k+\frac{\sigma(2c_0\sqrt{\log(nL)}+c_1\alpha_kt_n)}{h_k\sqrt
n},& k<k(n),\\
\displaystyle C_g\sqrt C_V\,(\rho_k+\delta_{k-1})^2h_k+\frac{\sigma(z
c_0+c_1\alpha_kt_n)}{h_k\sqrt n},& k=k(n),
\end{cases}
\\
\zeta_k&= 2\mu^*(\gamma_k^2\rho_k^{-2}+\sqrt2\,\gamma_k\rho_k^{-1}C_g),\\
\delta_k&= 2\gamma_k\sqrt{\mu^*}/\sqrt{1-\zeta_k},\\
\Omega_k&= \{\max_{\ell}|P_k^*(\hat\beta_{k,\ell}-\beta_{\ell})|\leq
\gamma_k\}.
\end{align*}
Here $\hat\beta_{k,\ell}=\frac1n\sum_{i=1}^n \widehat{\nabla
f}^{(k)}(X_i)\psi_\ell(X_i)$ with
$$
\begin{pmatrix}
\hat f^{(k)}(X_{i})\\[3pt]
\widehat{\nablaf}^{(k)}(X_{i})
\end{pmatrix}
=\bigg(\sum_{j=1}^{n}\onev{X_{ij}}\onev{X_{ij}}^{\T}\,w_{ij}^{(k)}\bigg)^{-1}
\sum_{j=1}^{n} Y_{j}\onev{X_{ij}}w_{ij}^{(k)},
$$
and $
w_{ij}^{(k)}=K\bigl(h_k^{-2}|(I+\rho_k^{-2}\widehat\Pi_{k-1})^{1/2}X_{ij}|^{2}\bigr)
$.

Combining Lemmas~\ref{lem5.3} and~\ref{lem5.4}, we obtain
$\Pb(\tr(I-\widehat\Pi_{k-1})\Pi^*>\delta_{k-1}^2)\leq
\Pb(\Omega_{k-1}^c)$ and therefore, using Corollary~\ref{cor1}, we get
\begin{align*}
\Pb\big(\Omega_k^c\big)&\le
\Pb\Big(\max_{\ell}|P_k^*(\hat\beta_{k,\ell}-\beta_{\ell})|>\gamma_k,\;
\tr(I-\widehat\Pi_{k-1})\Pi^*\leq\delta_{k-1}^2\Big)+\Pb\big(\Omega_{k-1}^c\big)\\
&\le \Pb\Big(\sup_{\Pi\in\mathscr
P_{m^*,\delta_{k-1}}}\max_{\ell}|P_k^*(\hat\beta_{k,\ell}-\beta_{\ell})|>\gamma_k\Big)+
\Pb\big(\Omega_{k-1}^c\big)\\
&\le  \frac3n+\Pb\big(\Omega_{k-1}^c\big),\qquad k\leq k(n)-1.
\end{align*}
Since $\Pb(\Omega_1^c)\leq 1/n$, it holds $\Pb(\Omega_{k(n)-1}^c)\le
(3k(n)-5)/n$ and
$\Pb(\Omega_{k(n)}^c)\le Lze^{-(z^2-1)/2}+\frac{3k(n)-5}{n}$.
Lemma~\ref{lem5.4} implies that
$$
\Pb\big(\tr(I-\widehat\Pi_{k(n)})\Pi^*>\delta_{k(n)}^2\big)\le
Lze^{-(z^2-1)/2}+\frac{3k(n)-5}{n}\ .
$$
According to Lemma~\ref{lem5.3}, we have
$\delta_{k(n)-2}\leq \rho_{k(n)-1}$, $\alpha_{k(n)-1}\leq 4$ and
$\zeta_{k(n)-1}\leq 1/2$. Consequently, for $n$ sufficiently large, we have
$$
\delta_{k(n)-1}=\frac{2\sqrt{\mu^*}\gamma_{k(n)-1}}
{\sqrt{1-\zeta_{k(n)-1}}}\leq C \bigg(\frac{\log(Ln)}{n}\bigg)^{1/2}\vee
n^{-2/3\vee m^*}
$$
and $\alpha_{k(n)}\leq 4\delta_{k(n)-1}\rho_{k(n)}^{-1}\leq
C[(\sqrt{\log(Ln)}(\rho_{k(n)}\sqrt{n})^{-1})\vee n^{-1/3\vee m^*}]$.
Since $h_{k(n)}=1$ and $(n\rho_{k(n)})^{-1}\le \rho_{k(n)}^2=n^{-2/(3\vee
m^*)}$, we
infer that
\begin{align*}
\gamma_{k(n)}&= C_g\sqrt C_V\,(\rho_{k(n)}+\delta_{k(n)-1})^2+\frac{\sigma(z
c_0+c_1\alpha_{k(n)}t_n)}{\sqrt n}\\
&\leq  Ct_n^2n^{-2/(3\vee m^*)}+\frac{c_0\sigma\, z}{\sqrt{n}}.
\end{align*}
Therefore
$\zeta_n:=\zeta_{k(n)}=O(\gamma_{k(n)}\rho_{k(n)}^{-1})$
tends to zero as
$n\to\infty$ at least as fast as $\sqrt{\log(nL)}\,n^{-1/(6\vee m^*)}$
and the assertion of the theorem follows from the definition of
$\delta_{k(n)}$ and Lemma~\ref{le4.2} (see below).

\subsection{Maximal inequality}\label{ss5}

The following result contains a well known maximal inequality for
the maximum of a Gaussian process. We include its proof for
the completeness of exposition.
\begin{proposition}\label{prop1}
Let $r$ be a positive number and let $\Gamma$ be a finite set. Let functions
$a_{j,\gamma}:\RR^p\to\RR^d$
obey the conditions
\begin{align}
&\sup_{\gamma\in\Gamma} \sup_{|u-u^*|\leq
r}\sum_{j=1}^n|a_{j,\gamma}(u)|^2\leq \kappa_0^2,\\
&\sup_{\gamma\in\Gamma} \sup_{|u-u^*|\leq r}\sup_{\e\in S_{d-1}}
\sum_{j=1}^n\bigg|\frac{d\
}{du}\;(\e^{\T}a_{j,\gamma}(u))\bigg|^2\leq\kappa_1^2.
\end{align}
If the $\varepsilon_j$'s are independent $\mathcal
N(0,\sigma^2)$-distributed random variables, then
$$
\Pb\bigg(\sup_{\gamma\in\Gamma}\sup_{|u-u^*|\leq r}\;\bigg|\sum_{j=1}^n
a_{j,\gamma}(u)\,\varepsilon_j\bigg|>t\sigma\kappa_0+2\sqrt
n\sigma\kappa_1\epsilon\bigg)\leq \frac2n,
$$
where $t=\sqrt{3\log(|\Gamma|(2r/\epsilon)^pn)}$.
\end{proposition}

\begin{Proof}
Let $B_r$ be the ball $\{u:|u-u^*|\leq r\}\subset\RR^p$ and
$\Sigma_{r,\epsilon}$ be the $\epsilon$-net on $B_r$
such that for any $u\in B_r$ there is an element $u_l\in\Sigma_{r,\epsilon}$
such that $|u-u_l|\leq\epsilon$. It
is easy to see that such a net with cardinality
$N_{r,\epsilon}<(2r/\epsilon)^p$ can be constructed. For every $u\in B_r$
we denote$\eta_\gamma(u)=\sum_{j=1}^n a_{j,\gamma}(u)\,\varepsilon_j$.
Since $\esp(|\eta_\gamma(u)|^2)\leq \sigma^2\kappa_0^2$ for any $\gamma$ and
for any $u$, we have
$$
\Pb\big(|\eta_\gamma(u_l)|>t\sigma\kappa_0\big)\leq
\Pb\Big(|\eta_\gamma(u_l)|>t\sqrt{\esp(|\eta_\gamma(u_l)|^2)}\Big)
\leq te^{-(t^2-1)/2}.
$$
Thus we get
\begin{align*}
\Pb\Big(\sup_{\gamma\in\Gamma}\sup_{u_l\in\Sigma_{r,\epsilon}}\;\big|\eta_\gamma(u_l)\big|>t\sigma\kappa_0\Big)
&\le \sum_{\gamma\in\Gamma}\sum_{l=1}^{N_{r,\epsilon}}
\Pb\Big(\big|\eta_\gamma(u_l)\big|>t\sigma\kappa_0\Big)\\
&\leq |\Gamma| N_{r,\epsilon}te^{-(t^2-1)/2}.
\end{align*}
Hence, if $t=\sqrt{3\log(|\Gamma|N_{r,\epsilon}n)}$, then
$
\Pb\Big(\sup_{\gamma\in\Gamma}\sup_{u_l\in\Sigma_{r,\epsilon}}\;\big|\eta_\gamma(u_l)\big|>t\sigma\kappa_0\Big)
\leq 1/n.
$
On the other hand, for any $u,u'\in B_r$ the Cauchy-Schwarz inequality
yields
\begin{align*}
\big|\eta_\gamma(u)-\eta_\gamma(u')\big|^2&= \sup_{\e\in S_{d-1}}
\big|\e^{\T}\big(\eta_\gamma(u)-\eta_\gamma(u')\big)\big|^2\\
&\le |u-u'|^2\cdot\sup_{u\in B_r}\sup_{\e\in
S_{d-1}}\bigg|\frac{d(\e^{\T}\eta_\gamma)}{du}\;(u)\bigg|^2\\
&=  |u-u'|^2\cdot\sup_{u\in B_r}\sup_{\e\in
S_{d-1}}\bigg|\sum_{j=1}^n\frac{d(\e^{\T}a_{j,\gamma})}{du}\;(u)
\,\varepsilon_j\bigg|^2\\
&\le  |u-u'|^2\cdot\sup_{u\in B_r}\sup_{\e\in
S_{d-1}}\sum_{j=1}^n\bigg|\frac{d(\e^{\T}a_{j,\gamma})}{du}\;(u)\bigg|^2
\sum_{j=1}^n\varepsilon_j^2\\
&\leq \kappa_1^2|u-u'|^2\sum_{j=1}^n\varepsilon_j^2.
\end{align*}
Since $\Pb\big(\sum_{j=1}^n \varepsilon_j^2>4n\sigma^2\big)$ is certainly
less than $n^{-1}$, we have
\begin{align*}
\Pb\Big(&\displaystyle\sup_{\gamma\in\Gamma}\sup_{u\in
B_r}\;\big|\eta_\gamma(u)\big|>t\sigma\kappa_0+2\sqrt
n\sigma\kappa_1\epsilon\Big)\\
&\le
\Pb\Big(\sup_{\gamma\in\Gamma}\sup_{u_l\in\Sigma_{r,\epsilon}}\;
\frac{|\eta_\gamma(u_l)|}{t\sigma\kappa_0}>1\Big)
+\Pb\Big(\sup_{\gamma\in\Gamma}\sup_{u\in
B_{r}}\;\frac{|\eta_\gamma(u)-\eta_\gamma(u_l(u))|}{2\sqrt
n\sigma\kappa_1\epsilon}>1\Big)\\
&\le
\frac1n
+\Pb\Big(\sup_{u\in
B_{r}}\kappa_1^2|u-u_l(u)|^2\sum_{j=1}^n\varepsilon_j^2>4
n\sigma^2\kappa_1^2\epsilon^2\Big)\leq
\frac2n,
\end{align*}
and the assertion of proposition follows.
\end{Proof}

\subsection{Properties of the solution to (\ref{maxopt})}\label{ss6}

We collect below some simple facts concerning the solution to
the optimization problem (\ref{maxopt}). By classical arguments, it
is always possible to choose a measurable solution $\widehat\Pi$
to (\ref{maxopt}). This measurability will be assumed in the sequel.

In Proposition~\ref{P1} the case of general $m$ (not necessarily
equal to $m^*$) is considered. As we explain below, this generality
is useful for further developments of the method extending it to the
case of unknown structural dimension $m^*$.

The vectors $\beta_\ell$ are assumed to belong to a $m^*$-dimensional
subspace $\cc{S}$ of $\RR^d$, but in this subsection we do not assume
that $\beta_\ell$s are defined by (\ref{betal}). In fact, we will
apply the results of this subsection to the vectors
$\Pi^*\hat\beta_\ell$.

Denote
\begin{align*}
R(\Pi)&=\max_{\ell} {\hat\beta}_{\ell}^{\T} (I-\Pi) \hat\beta_{\ell},\\
\hat{\cc{R}}(m) &= \min_{\Pi \in \cc{A}_{m}} \sqrt{R(\Pi)} =
\sqrt{R(\widehat{\Pi}_{m})}.
\end{align*}
We also define
\begin{align*}
\cc{R}^*(m)&=\min_{\Pi \in \cc{A}_{m}} \max_\ell|(I-\Pi)^{1/2}\beta_\ell|.
\end{align*}
and denote by $\Pi_m^*$ a minimizer of
$\max_\ell\beta_\ell^\T(I-\Pi)\beta_\ell$ over $\Pi\in\cc{A}_m$.
Since for $m\ge m^*$ the projector $\Pi^*$ is in $\cc{A}_m$, we have
$\Pi_m^*=\Pi^*$ and $\cc{R}^*(m)=0$.

\begin{proposition}
\label{P1}
Let $\B =\bigl\{\betam = \sum_{\ell}c_{\ell}\beta_\ell : \sum_{\ell}
|c_{\ell}|
\le 1 \bigr\}$ be the convex hull of vectors $\beta_\ell$.
If $\max_\ell |\hat\beta_\ell-\beta_\ell|\le \eps$, then
\begin{align*}
\hat{\cc{R}}(m) &\le  \cc{R}^*(m)+\varepsilon,\\
\max_{\betam\in\B} |(I-\widehat{\Pi}_{m})^{1/2} \betam|&\le
\cc{R}^*(m)+2\varepsilon.
\end{align*}
When $m< m^{*}$, we have also the lower bound
$\hat{\cc{R}}(m)\ge (\cc{R}^{*}(m) - \varepsilon)_{+}$.
\end{proposition}
\begin{Proof}
For every $\ell\in{1,\ldots,L}$, we have
\begin{align*}
|(I-\Pi^{*}_{m})^{1/2} \hat\beta_{\ell}|
&\le |(I-\Pi^{*}_{m})^{1/2} \beta_\ell|
+ |(I-\Pi^{*}_{m})^{1/2}(\hat\beta_{\ell} - \beta_\ell)|\\
&\le \cc{R}^{*}(m) + |\hat\beta_{\ell} - \beta_\ell|
\le \cc{R}^{*}(m) + \varepsilon.
\end{align*}
Since $\widehat\Pi_m$ minimizes $\max_\ell |(I-\Pi)^{1/2}\hat\beta_\ell|$
over
$\Pi\in\mathcal A_m$, we have
$$
\max_\ell |(I-\widehat\Pi_m)^{1/2}\hat\beta_\ell|\leq \max_\ell
|(I-\Pi^*_m)^{1/2}\hat\beta_\ell|\leq
\cc{R}^{*}(m)+ \varepsilon.
$$
Denote $ A = (I - \widehat{\Pi}_{m})^{1/2}$. From definition $ 0 \preceq A
\preceq I $.
Therefore, for every $\ell $
\begin{align*}
| A \beta_\ell|\le|A \hat\beta_{\ell}| + |A (\beta_\ell - \hat\beta_{\ell})|
\le|A \hat\beta_{\ell}| + |\beta_\ell - \hat\beta_{\ell}|
\le\cc{R}^{*}(m)+ 2 \varepsilon.
\end{align*}
The second inequality of the proposition follows now from
$|A \betam|\leq \max_\ell| A \beta_\ell|$ for every $\betam\in\B$.

To prove the last assertion, remark that according to the
definition of $\cc{R}^{*}(m)$, for every matrix $\Pi \in
\cc{A}_{m}$ there exists an index $\ell $ such that
$ |(I - \Pi)^{1/2} \beta_\ell| \ge \cc{R}^{*}(m)$. In
particular, $ |(I - \widehat{\Pi}_{m})^{1/2} \beta_\ell| \ge
\cc{R}^{*}(m)$ for some $\ell$ and hence
$|(I - \widehat{\Pi}_{m})^{1/2} \hat\beta_{\ell}|
  \ge |(I - \widehat{\Pi}_{m})^{1/2} \beta_\ell| - |\hat\beta_{\ell} -
\beta_\ell|\ge \cc{R}^{*}(m) - \varepsilon$.
\end{Proof}

Proposition~\ref{P1} can be used for estimating the
structural dimension $m$.
Indeed, $\hat{\cc{R}}(m) \le \varepsilon$ for $ m \ge m^{*}$ and
the results mean that $\hat{\cc{R}}(m)\ge(\cc{R}^{*}(m)-\varepsilon)_{+}$
for $m< m^{*}$. Therefore, it is natural to search for the smallest
value $\widehat{m}$ of $ m $ such that the function $\hat{\cc{R}}(m)$
does not significantly decrease for $m\ge \widehat{m}$.

From now on, we assume that the structural dimension $m^*$ is known
and write $\widehat\Pi$ instead of $\widehat\Pi_{m^*}$.

\begin{proposition}
\label{P2}
If the vectors $\beta_\ell$ satisfy \textbf{(A2)} and
$\max_\ell|\hat\beta_\ell-\beta_\ell|\le\varepsilon$, then
$\tr (I - \widehat{\Pi}) \Pi^{*}\le 4 \varepsilon^{2}  \mu^*$
and $\tr[ (\widehat{\Pi} - \Pi^{*})^{2}]\le 8 \varepsilon^{2}\mu^{*}$.
\end{proposition}

\begin{Proof}
In view of the relations
$\tr \widehat{\Pi}^{2} \le \tr \widehat{\Pi}\le  m^*$ and
$\tr (\Pi^{*})^{2} = \tr \Pi^{*} = m^{*}$, we have
\begin{align*}
\tr (\widehat{\Pi} - \Pi^{*})^{2}&=
\tr (\widehat{\Pi}^{2} - \Pi^{*}) + 2 \tr (I - \widehat{\Pi}) \Pi^{*}
\le 2 | \tr (I-\widehat{\Pi}) \Pi^{*}|.
\end{align*}
Note also that the equality
$\tr (I-\widehat{\Pi})
\Pi^{*}=\tr(I-\widehat\Pi)^{1/2}\Pi^*(I-\widehat\Pi)^{1/2}$
implies that $\tr (I-\widehat{\Pi}) \Pi^{*}\geq 0$.
Now condition (\ref{Pi*}) and Proposition~\ref{P1} imply
\begin{align*}
\tr (I-\widehat{\Pi})
\Pi^{*}&=\tr(I-\widehat\Pi)^{1/2}\Pi^*(I-\widehat\Pi)^{1/2}\\
& \le
\sum_{k=1}^{m^{*}} \mu_{k} \tr (I-\widehat{\Pi})^{1/2} \betam_{k}
\betam_{k}^{\T}(I-\widehat{\Pi})^{1/2}\\
&\le
\sum_{k=1}^{m^{*}} \mu_{k} \betam_{k}^{\T} (I-\widehat{\Pi}) \betam_{k}
\le
(2 \varepsilon)^{2} \sum_{k=1}^{m^{*}} \mu_{k}
\end{align*}
and the assertion follows.
\end{Proof}

\begin{lemma}
\label{LAB}
Let $\tr (I - \widehat{\Pi}) \Pi^{*} \le \delta^{2}$
for some $\delta < 1 $. Then
for any $ x \in \RR^{d}$
\begin{align*}
|\Pi^{*} x| \le |\widehat{\Pi}^{1/2} x| + \delta |x| .
\end{align*}
\end{lemma}
\begin{Proof}
Denote $\hat{A} = \widehat{\Pi}^{1/2}$.
It obviously holds
$ |\Pi^{*} x| \le |\Pi^{*} \hat{A} x| + |\Pi^{*} (I - \hat{A}) x|
$ and
\begin{align*}
|\Pi^{*} (I - \hat{A}) x|^{2}
\le \|\Pi^*(I-\hat A)\|_2^2\cdot|x|^{2} \le
\tr[ \Pi^{*} (I - \hat{A})^{2} \Pi^{*}]\cdot |x|^{2}.
\end{align*}
For every $\Pi \in \cc{A}_{m}$, it obviously holds
$ (I - \Pi^{1/2})^{2} = I - 2 \Pi^{1/2} + \Pi \preceq I - \Pi $, and
hence,
$\tr \Pi^{*} (I - \Pi^{1/2})^{2} \Pi^{*} \le \tr \Pi^{*} (I - \Pi) \Pi^{*}
$.
Therefore,
\begin{align*}
\tr \Pi^{*} (I - \hat{A})^{2} \Pi^{*}
\le
\tr \Pi^{*} (I - \widehat{\Pi}) \Pi^{*}
= \tr (I - \widehat{\Pi}) \Pi^{*} \le  \delta^{2}
\end{align*}
yielding
$|\Pi^{*} x| \le |\Pi^{*} \hat{A} x| + \delta |x|
\le |\hat{A} x| + \delta |x|$ as required.
\end{Proof}

\begin{corollary}
\label{Cellipt}
Let $\rho \in (0,1) $,
and
$\hat{S}_{\rho}  = (I + \rho^{-2} \widehat{\Pi})^{1/2}$.
If $\tr (I - \widehat{\Pi}) \Pi^{*} \le \delta^{2}$,
then for any $ x \in \RR^{d}$, the condition $ |\hat{S}_{\rho} x| \le h
$ implies
$ |\Pi^{*} x| \le (\rho + \delta) h $.
\end{corollary}

\begin{Proof}
The result follows from Lemma~\ref{LAB} and the obvious inequalities
$ |x| \le |\hat{S}_{\rho} x| \le h $ and
$ |\widehat{\Pi}^{1/2} x| \le \rho|\hat{S}_{\rho} x| \le \rho h $.
\end{Proof}

\begin{lemma}\label{le4.2}
Let $\tr (I - \widehat{\Pi}) \Pi^{*} \le \delta^{2}$
for some $\delta\in[0,1[$ and let $\widehat\Pi_{m^*}$ be the orthogonal
projection matrix in $\RR^d$ onto the subspace spanned by the eigenvectors
of $\widehat\Pi$ corresponding to its largest $m^*$ eigenvalues.
Then $\tr (I - \widehat{\Pi}_{m^*}) \Pi^{*} \le \delta^{2}/(1-\delta^2)$.
\end{lemma}
\begin{Proof}
Let $\hat\lambda_j$ and $\hat\theta_j$, $j=1,\ldots,d$ be respectively
the eigenvalues and the eigenvectors of $\widehat\Pi$. Assume that
$\hat\lambda_1\ge \hat\lambda_2\ge \ldots\ge \hat\lambda_d$. Then
$\widehat\Pi=\sum_{j=1}^d\hat\lambda_j \hat\theta_j\hat\theta_j^\T$ and
$\widehat\Pi_{m^*}=\sum_{j=1}^{m^*}\hat\theta_j\hat\theta_j^\T$. Moreover,
$\sum_{j=1}^d \hat\theta_j\hat\theta_j^\T=I$ since
$\{\hat\theta_1,\ldots,\hat\theta_d\}$ is an orthonormal basis of $\RR^d$,
Therefore,
on the one hand,
\begin{align*}
\tr[\widehat\Pi\Pi^*]&\le \sum_{j\le m^*} \hat\lambda_j
\tr[\hat\theta_j\hat\theta_j^\T \Pi^*]+
\hat\lambda_{m^*}\sum_{j> m^*}\tr[\hat\theta_j\hat\theta_j^\T \Pi^*]\\
&=\sum_{j\le m^*} (\hat\lambda_j-\hat\lambda_{m^*})\tr[\hat\theta_j\hat\theta_j^\T \Pi^*]+
\hat\lambda_{m^*}\tr\Big[\sum_{j=1}^d\hat\theta_j\hat\theta_j^\T \Pi^*\Big]\\
&=\sum_{j\le m^*} (\hat\lambda_j-\hat\lambda_{m^*})
\tr[\hat\theta_j\hat\theta_j^\T \Pi^*]+m^*\hat\lambda_{m^*}.
\end{align*}
Since $\tr[\hat\theta_j\hat\theta_j^\T \Pi^*]=|\Pi^*\hat\theta_j|^2\le 1$, we
get $\tr[\widehat\Pi\Pi^*]\le \sum_{j\le m^*}\hat\lambda_{j}$. Taking into
account the relations $\sum_{j\le d} \hat\lambda_j\le m^*$, $\tr\Pi^*=m^*$
and $(1-\hat\lambda_{m^*+1})(I-\widehat\Pi_{m^*})\preceq I-\widehat\Pi$,
we get $\lambda_{m^*+1}\le m^*-\sum_{j\le m^*}\hat\lambda_j \le
\tr[(I-\widehat\Pi)\Pi^*]\le\delta^2$ and therefore
$\tr[(I-\widehat\Pi_{m^*})\Pi^*]\le \delta^2/(1-\hat\lambda_{m^*+1})\le
\delta^2/(1-\delta^2)$.
\end{Proof}

\subsection{Technical lemmas}\label{ss7}

This subsection contains five technical results. The first three lemmas
have been used in the proof of Proposition~\ref{P5}, whereas the two
last lemmas have been used in the proof of Theorem~\ref{THM1}.

\begin{lemma}\label{lemalpha}
If $\rho\leq 1$, then $\|U-U^*\|_2\leq \alpha$.
\end{lemma}
\begin{Proof} The inequality $P_\rho^*\preceq (I-\Pi^*)+\rho \Pi^*$ implies
that
\begin{align*}
\rho^{2}\big\|U-U^*\big\|_2&=
\big\|P_{\rho}^{*}(\Pi-\Pi^*)P_{\rho}^{*}\big\|_2 \\
&\le
\rho^2\big\|\Pi^*(\Pi-\Pi^*)\Pi^*\big\|_2+\big\|(I-\Pi^*)(\Pi-\Pi^*)(I-\Pi^*)\big\|_2\\
&+2\rho\big\|\Pi^*(\Pi-\Pi^*)(I-\Pi^*)\big\|_2.
\end{align*}
Since $\|A\|_2^2=\tr AA^\T \leq (\tr (AA^\T)^{1/2})^2$ for any matrix $A$,
it holds
\begin{align*}
\big\|\Pi^*(\Pi-\Pi^*)\Pi^*\big\|_2&= \big\|\Pi^*(I-\Pi)\Pi^*\big\|_2\\
&\le \tr\, \Pi^*(I-\Pi)\Pi^*=\tr (I-\Pi)\Pi^*\leq \delta^2.
\end{align*}
By similar arguments one checks that
\begin{align*}
\big\|(I-\Pi^*)(\Pi-\Pi^*)(I-\Pi^*)\big\|_2
&= \big\|(I-\Pi^*)\Pi(I-\Pi^*)\big\|_2\leq \tr (I-\Pi^*)\Pi\\
&= \tr\Pi-\tr\Pi^*+\tr\Pi^*(I-\Pi)\\
&\le  m^*-m^*+\delta^2,\\
\big\|\Pi^*(\Pi-\Pi^*)(I-\Pi^*)\big\|_2&\le
\big\|\Pi^*(\Pi-\Pi^*)\big\|_2=\big\|\Pi^*(I-\Pi)\big\|_2\\
&\le \big\|\Pi^*(I-\Pi)^{1/2}\big\|_2\leq(\tr\Pi^*(I-\Pi)\Pi^*)^{1/2}\\
&=  (\tr(I-\Pi)\Pi^*)^{1/2}\leq
\delta.
\end{align*}
Thus we get $\big\|U-U^*\big\|_2\leq
\delta^2(1+\rho^{-2})+2\delta\rho^{-1}$.
The assumption $\rho\leq 1$ yields the assertion of the lemma.
\end{Proof}

\begin{lemma}\label{lem4.1}
If $\psi_\ell$s and $U$ satisfy \textbf{(A3)} and (\ref{asspsi}), then
$$
\sum_{j=1}^n|c_{j,\ell}(U)|^2\leq \frac{dC_KC_V\bar\psi^2}{h^2n}.
$$
\end{lemma}
\begin{Proof}
Simple computations yield
\begin{align}
\sum_{j=1}^n\bigg|\V_{i}^{-1} \onev{Z_{ij}} \bigg|^2w_{ij}&
=\tr(\V_i^{-1})\leq \frac{d C_V}{N_i}.\qquad\label{dCV}
\end{align}
Hence, we have
\begin{align*}
\sum_{j=1}^n|c_{j,\ell}|^2&= \frac{1}{h^2n^2}\sum_{j=1}^n
\bigg|\sum_{i=1}^{n} \V_{i}^{-1} \onev{Z_{ij}}
w_{ij}\,\psi_{\ell}(X_{i})\bigg|^2\\
&\le \frac{\bar\psi^2}{h^2n^2}\:
\sum_{j=1}^n\bigg(\sum_{i=1}^n\frac{w_{ij}}{N_i}\bigg)
\bigg(\sum_{i=1}^n\bigg|\V_{i}^{-1} \onev{Z_{ij}}\bigg|^2 N_iw_{ij}\bigg)\\
&\le  \frac{C_K\bar\psi^2}{h^2n^2}\:
\sum_{j=1}^n\sum_{i=1}^n\bigg|\V_{i}^{-1} \onev{Z_{ij}}\bigg|^2 N_iw_{ij}.
\end{align*}
Interchanging the order of summation and using inequality (\ref{dCV}) we get
the desired result.
\end{Proof}

\begin{lemma}\label{lem4.2}
If \textbf{(A3)} and (\ref{asspsi})
are fulfilled, then, for any $j=1,\ldots,n$,
\begin{equation*}
\sup_U\sup_{\e\in S_{d-1}}\bigg|\frac{d\
}{dU}\,(\e^\T c_{j,\ell})(U)\bigg|^2\leq
C\bigg(\frac{C_w^2C_V^4C_K^2\bar\psi^2}{n^2h^2}
+\frac{C_V^2C_{K'}^2\bar\psi^2}{n^2h^2}\bigg),
\end{equation*}
where $C$ is a numerical constant and $\frac{d\,}{dU}(\e^{\T}
c_{j,\ell})(U)$ is the $d\times d$ matrix with entries
$\frac{\partial\,\e^{\T} c_{j,\ell}(U)}{\partial U_{pq}}$.
\end{lemma}
\begin{Proof}
We have
\begin{align*}
\bigg\|\frac{d\,\e^{\T} c_{j,\ell}(U)}{dU}\bigg\|_2^2&\le
2\bigg\|\frac{1}{hn}\sum_{i=1}^{n}\bigg[\frac{d}{dU}\;\e^{\T}\V_{i}^{-1}(U)
\onev{Z_{ij}}\bigg]
w_{ij}(U)\psi_{\ell}(X_{i})\bigg\|_2^2\\
&+
2\bigg\|\frac{1}{hn}\sum_{i=1}^{n}\e^{\T}\V_{i}^{-1}(U) \onev{Z_{ij}}\frac{d
w_{ij}(U)}{dU}\;\psi_{\ell}(X_{i})
\bigg\|_2^2\\
&= \Delta_1+\Delta_2.
\end{align*}
One checks that $\|d\,w_{ij}(U)/dU\|_2=|w_{ij}'(U)|\cdot |Z_{ij}|^2\leq
5|w_{ij}'(U)|$, where we used the notation
$w_{ij}'(U)=K'(Z_{ij}^{\T}U Z_{ij})$ and the inequality
\begin{align*}
h^2|Z_{ij}|^2&=|S_\rho^*X_{ij}|^2=|(I-\Pi^*)X_{ij}|^2+2\rho^{-2}|\Pi^*X_{ij}|^2\\
&\leq h^2+2(\delta/\rho+1)^2h^2\leq 5h^2,
\end{align*}
which follows from Lemma~\ref{LAB}. We get
\begin{align*}
\Delta_2
&\le
\frac{50\bar\psi^2}{n^2h^2}\bigg(\sum_{i=1}^{n}\Big|\V_{i}^{-1}(U)\onev{Z_{ij}}w'_{ij}(U)\Big|
\bigg)^2\leq \frac{C \bar\psi^2 C_V^2 C_{K'}^2}{n^2h^2}.
\end{align*}
In order to estimate the term $\Delta_1$, remark that the differentiation
(with respect to $U_{pq}$)
of the identity $\V_i^{-1}(U) \V_i(U)=I_{d+1}$ yields
$$
\frac{\partial \V_i^{-1}}{\partial U_{pq}}\;(U)=-\V_i^{-1}(U)\;
\frac{\partial \V_i}{\partial U_{pq}}\;(U)\V_i^{-1}(U).
$$
Simple computations show that
\begin{align*}
\frac{\partial \V_i}{\partial U_{pq}}\;(U)&= \sum_{j=1}^{n}
\onev{Z_{ij}}
\onev{Z_{ij}}^{\T}
\frac{\partial \ }{\partial U_{pq}}w_{ij}(U)\\
&= \sum_{j=1}^{n}
\onev{Z_{ij}}
\onev{Z_{ij}}^{\T}w'_{ij}(U)
(Z_{ij})_p(Z_{ij})_q.
\end{align*}
Hence, for any $a_1,a_2\in\RR^{d+1}$,
\begin{align*}
\frac{d a_1^{\T}\V_i^{-1}a_2}{d
U}\;(U)&= \sum_{j=1}^{n}a_1^{\T}\V_i^{-1}(U)\onev{Z_{ij}}
\onev{Z_{ij}}^{\T}\V_i^{-1}(U)a_2\;w'_{ij}(U)Z_{ij}Z_{ij}^{\T}.
\end{align*}
This relation combined with the estimate $|Z_{ij}|\leq 5$ for all $i,j$ such
that
$w_{ij}\not=0$, implies the norm estimate
\begin{align*}
\bigg\|\frac{d a_1^{\T}\V_i^{-1}a_2}{d U}\;(U)\bigg\|_2&\leq
25\sum_{j=1}^{n}
\bigg|a_1^{\T}\V_i^{-1}(U)\onev{Z_{ij}}\onev{Z_{ij}}^{\T}\V_i^{-1}(U)a_2\;
w'_{ij}(U)\bigg|\\
&\leq  150|a_1|\,|a_2|\sum_{j=1}^{n}
\big\|\V_i^{-1}(U)\big\|^2|w'_{ij}(U)|\\
&\le  150 C_w C_V^2|a_1|\,|a_2| N_i(U)^{-1}.
\end{align*}
It leads to the estimate $\Delta_1\leq C\,\frac{C_w^2 C_V^4 C_K^2\bar\psi^2}{n^2h^2}$,
and the assertion of the lemma follows.
\end{Proof}

\begin{lemma}\label{lem5.3}
There exists an integer $n_0\ge 0$ such that, as soon as $n\ge n_0$,
$\delta_{k-1}\le \rho_k$, $\alpha_k\leq 4$ and $\zeta_k\leq 1/2$ for
all $k\in\{2,\ldots,k(n)\}$.
\end{lemma}
\begin{Proof}
In view of $C_0n^{-1/(d\vee 4)}=\rho_1 h_1$ and
$\rho_{k(n)}h_{k(n)}\ge C_2 n^{-1/3}$, the sequence
$$
s_n=4\sqrt{C_V} C_g h_1+
\frac{4\sigma(c_0\sqrt{\log(Ln)}+c_1t_n)}{\sqrt{n}\,\rho_{k(n)}h_{k(n)}}
$$
tends to zero as $n\to\infty$.

We do now an induction on $k$.
Since $s_n\to 0$ as $n\to\infty$ and $\gamma_1\le s_n$,
the inequality $\delta_1=2\gamma_1\sqrt{\mu^*}\leq
1/\sqrt{2}=\rho_1/\sqrt{2}$ is true for sufficiently large
values of $n$.
Let us prove the implication
$$
\delta_{k-1}\leq \rho_{k-1}/\sqrt 2\quad\Longrightarrow\quad
\begin{cases}
\zeta_k\leq 1/2,\\
\delta_k\leq \rho_k/\sqrt 2.
\end{cases}
$$
Since $1/\sqrt2\le e^{-1/6}$ we infer that $\delta_{k-1}\le \rho_k$ and
therefore $\alpha_k\le 4$.
By our choice of $a_h$ and $a_\rho$, we have
$\rho_1 h_1\ge \rho_k h_k\ge \rho_{k(n)}h_{k(n)}$.
Therefore,
\begin{align*}
\frac{\gamma_k}{\rho_k}&\le  4\sqrt{C_V} C_g\rho_k
h_k+\frac{4\sigma(c_0\sqrt{\log(Ln)}+c_1t_n)}{\sqrt n\,\rho_kh_k}\\
&\le
4\sqrt{C_V} C_gh_1+\frac{4\sigma(c_0\sqrt{\log(Ln)}+c_1t_n)}
{\sqrt{n}\rho_{k(n)}h_{k(n)}}=s_n.
\end{align*}
Thus, for $n$ large enough,
$\zeta_k\leq 1/2$ and $\gamma_k\leq \rho_k/4$. This implies that
$\delta_k=2\gamma_k(1-\zeta_k)^{-1/2}\leq \rho_k/\sqrt 2$.
By induction we infer that $\delta_{k-1}\leq \rho_{k-1}/\sqrt2\leq \rho_{k}$
and $\zeta_k\leq 1/2$ for any
$k=2,\ldots,k(n)-1$. This completes the proof of the lemma.
\end{Proof}

\begin{lemma}\label{lem5.4}
If $k>2$ and $\zeta_{k-1}<1$ then
$\Omega_{k-1}\subset \{\tr (I-\widehat\Pi_{k-1})\Pi^*\leq \delta_{k-1}^2\}$.
\end{lemma}
\begin{Proof}
Let us denote $\tilde\beta_{\ell}=\Pi^*\hat\beta_{{k-1},\ell}$, then
$\tilde\beta_{\ell}\in\cc{S}^*$ and
under $\Omega_{k-1}$ we have
$$
|P_{k-1}^*(\hat\beta_{{k-1},\ell}-\beta_{\ell})|\leq
\gamma_{{k-1}}\quad\Longrightarrow\quad
\begin{cases}
\max_\ell |\hat\beta_{{k-1},\ell}-\tilde\beta_{\ell}|\leq
\gamma_{{k-1}},\\
\max_\ell |\tilde\beta_{\ell}-\beta_{\ell}|\leq \sqrt
2\gamma_{{k-1}}/\rho_{k-1}.
\end{cases}
$$
Set $B=\sum_{i=1}^{m^*}\mu_{i}\bar\beta_{i}\bar\beta_{i}^\T$
and $\tilde
B=\sum_{i=1}^{m^*}\mu_{i}\bar{\tilde\beta}_{i}\bar{\tilde\beta}_{i}^\T$,
where $\bar{\tilde\beta}_{i}=\sum_{\ell}
c_\ell\tilde\beta_{\ell}$ if
$\bar{\beta}_{i}=\sum_{\ell}c_\ell\beta_{\ell}$. Since $\sum_\ell
|c_\ell|\leq 1$, we have
$|\bar\beta_{i}|\leq \max_\ell |\beta_{\ell}|\leq\|\nabla f\|_\infty $ and
$|\bar\beta_{i}-\bar{\tilde\beta}_{i}|\leq
\max_\ell |\beta_{\ell}-\tilde\beta_{\ell}|$. Therefore
\begin{align*}
\|B-\tilde B\|&\le \sum_{{i}=1}^{m^*}
\mu_{i}\|\bar{\beta}_{i}\bar{\beta}_{i}^\T-\bar{\tilde\beta}_{i}\bar{\tilde\beta}_{i}^\T\|
\le
\mu^*\max_k\|\bar{\beta}_{i}\bar{\beta}_{i}^\T-\bar{\tilde\beta}_{i}\bar{\tilde\beta}_{i}^\T\|\\
&\le
\mu^*\max_{i}\Big(|\bar{\beta}_{i}-\bar{\tilde\beta}_{i}|^2+2|\bar{\beta}_{i}|\cdot|\bar{\beta}_{i}-\bar{\tilde\beta}_{i}|\Big)\\
&\le
\mu^*\big(2\gamma_{{k-1}}^2\rho_{k-1}^{-2}+2\sqrt{2}\,\gamma_{k-1}\rho_{k-1}^{-1}\max_\ell
|\beta_{\ell}|\big)=\zeta_{k-1}
\end{align*}
and hence, for every unit vector $v\in \cc{S}^*$,
$v^{\T}\tilde Bv\ge \big(v^{\T}Bv-\big|v^{\T}Bv-v^{\T}\tilde Bv\big|\big)
\ge v^{\T}Bv-\|B-\tilde B\|\ge 1-\zeta_{k-1}$.
This inequality implies that $\Pi^*\preceq (1-\zeta_{k-1})^{-1}\tilde B$
and, in view of Proposition~\ref{P2}
we obtain the assertion of the lemma.
\end{Proof}

\textbf{Acknowledgment.} Much of this work has been carried out when the first
author was visiting the Weierstrass Institute for Applied Analysis and
Stochastics. The financial support from the institute and the hospitality
of Professor Spokoiny are gratefully acknowledged.



\begin{thebibliography}{9}

\bibitem{ALL} {\sc  Antoniadis, A., Lambert-Lacroix, S. and Leblanc, F.} (2003). Effective
Dimension reduction methods for tumor classification using gene
expression data, \emph{Bioinformatics}, 19, 563--570.

\bibitem{Bell} {\sc Bellman, R.E.} \emph{Adaptive Control Processes.}
Princeton University Press, Princeton, NJ, 1961.

\bibitem{BKRW}  {\sc Bickel, P., Klaassen, C., Ritov, Y.\ and Wellner, J.}\
(1998).
\textit{Efficient and Adaptive Estimation for Semiparametric Models,%
} Springer, New York.

\bibitem{BuraCook} {\sc Bura, E.\ and Cook, R. D.} (2001). Estimating the
structural dimension of regressions via parametric inverse regression.
\emph{J. R. Stat. Soc. Ser. B Stat. Methodol.}  63  (2001),  no. 2,
393--410.

\bibitem{BuPf} \textsc{Bura, E. and Pfeiffer, R. M.} (2003), Graphical Methods for Class
Prediction Using Dimension Reduction Techniques on DNA Microarray
Data, \emph{Bioinformatics}, 19, 1252--1258.

\bibitem{ChLiT} {\sc Chan, K. S., Li, M. C. and Tong, H.} (2004).
Partially Linear Reduced-rank Regression . Technical report.  \texttt{www.stat.uiowa.edu/techrep/tr328.pdf}

\bibitem{Cook} {\sc Cook, R. D.} \emph{Regression graphics. Ideas for
studying
regressions through graphics.} Wiley Series in Probability and Statistics:
Probability and Statistics. John Wiley \& Sons, Inc., New York, 1998.

\bibitem{CookLi} {\sc Cook, R. D. and Li, B.} (2002). Dimension reduction
for conditional mean in regression.  \textit{Ann. Statist.} 30,
no. 2, 455--474.

\bibitem{CookLi04} {\sc Cook, R. D. and Li, B.} (2004). Determining the
dimension of iterative Hessian transformation.
\emph{Ann. Statist.}  32,  no. 6, 2501--2531.

\bibitem{CookNi} {\sc Cook, R. D. and Ni, L.} (2005). Sufficient dimension
reduction via
inverse regression: a minimum discrepancy approach.
\emph{J. Amer. Statist. Assoc.}  100,  no. 470, 410--428.

\bibitem{CookWeis1} {\sc Cook, R. D. and Weisberg, S.} (1991).
Discussion of ``Sliced inverse regression for dimension reduction''.by
K.\,C.\,Li, \emph{J. Amer. Statist. Assoc.}  86, no.\ 414, 328--332.

\bibitem{CookWeis} {\sc Cook, R. D. and Weisberg, S.} \emph{Applied Regression Including
Computing and Graphics.} Hoboken NJ: John Wiley, 1999.

\bibitem{DHP}
{\sc Delecroix, M., Hristache, M. and Patilea, V.} (2006). On semiparametric
$M$-estimation in
single-index regression.  .
\emph{J. Statist. Plann. Inference }  136,  no. 3, 730--769.

\bibitem{FanGij} {\sc Fan, J. and Gijbels, I.} \emph{Local polynomial
modelling
and its applications.} Monographs on Statistics and Applied Probability, 66.
Chapman \& Hall, London, 1996.

\bibitem{HJPS} {\sc Hristache, M., Juditsky, A.,
Polzehl, J. and Spokoiny, V.} (2001). Structure adaptive approach for
dimension reduction. \emph{Ann. Statist.}  29,  no. 6, 1537--1566.

\bibitem{HJS} {\sc Hristache, M., Juditsky, A. and Spokoiny, V.} (2001).
Direct estimation of the index coefficient in a single-index model.
\emph{Ann. Statist.}  29,  no. 3, 595--623.

\bibitem{LiDuan} {\sc Li, K. C. and Duan, N.} (1989). Regression analysis
under link
violation.  \textit{Ann. Statist.}  17,  no. 3, 1009--1052.


\bibitem{LiSIR} {\sc Li, K. C.} (1991). Sliced inverse regression for
dimension reduction. With discussion and a rejoinder by the author.
\emph{J. Amer. Statist. Assoc.}  86  (1991),  no. 414, 316--342.


\bibitem{Li} {\sc Li, K. C.} (1992).  On principal hessian directions for
data visualization and dimension reduction: another application of Stein's
lemma. \emph{J. Amer. Statist. Assoc.}, 87, 1025--1039.

\bibitem{Sam} {\sc Samarov, A.}   (1993). Exploring regression structure
using nonparametric functional estimation.
\emph{J. Amer. Statist. Assoc.}  88,  no. 423, 836--847.

\bibitem{SSV} {\sc Samarov, A., Spokoiny, V. and Vial, C.} (2006)
Component identification and estimation in nonlinear high-dimensional
regression models by structural adaptation.
\emph{J. Amer. Statist. Assoc.} 100 (2005), no. 470, 429--445.

\bibitem{XiaetAl} {\sc Xia, Y., Tong, H., Li, W. K. and Zhu, L. X.} (2002).
An adaptive estimation of dimension reduction space.  \emph{J. R. Stat. Soc.
Ser. B Stat. Methodol.}  64,  no. 3, 363--410.

\bibitem{YinCook} {\sc Yin, X. and Cook, R. D.} (2005).
Direction estimation in single-index regressions.
\textit{Biometrika} 92,  no. 2, 371--384.

\end{thebibliography}
\end{document}